\documentclass{amsart}
\usepackage[left=1.5in,right=1.5in,bottom=1.5in]{geometry}
\usepackage{amssymb}
\usepackage{amsthm}
\usepackage{amsmath}
\usepackage{hyperref}
\usepackage{comment}

\newtheorem{thm}{Theorem}[section]
\newtheorem{lem}[thm]{Lemma}

\newtheorem{corollary}[thm]{Corollary}
\newtheorem{prop}[thm]{Proposition}
\newtheorem{proposition}[thm]{Proposition}

\newtheorem{algorithm}[thm]{Algorithm}

\theoremstyle{definition}
\newtheorem{remark}[thm]{Remark}
\newtheorem{defn}[thm]{Definition}
\newtheorem{example}[thm]{Example}

\usepackage{todonotes}
\usepackage{mathtools}
\newcommand{\Red}{\mathcal{R}}
\newcommand{\Ans}{{\mathcal{M}}}
\newcommand{\AnsD}{{\mathcal{\bar{M}}}}
\newcommand{\Gleb}{arc }
\newcommand{\TitleGleb}{Arc }
\newcommand{\col}{\operatorname{col}}
\newcommand{\wt}{\operatorname{wt}}
\newcommand{\wtlim}{\operatorname{wtlim}}
\renewcommand{\tilde}{\widetilde}
\newcommand{\id}{\mathrm{id}}
\DeclareMathOperator{\code}{code}
\definecolor{pumpkin}{RGB}{255,117,24}

\begin{document}

\title[The maximum multiplicity of a generator in a reduced word]{Repeatable patterns and the maximum multiplicity of a generator in a reduced word}

\author[Gaetz]{Christian Gaetz} 
\thanks{C.G. was partially supported by an NSF Fellowship under grant DMS-2103121.}
\address{Department of Mathematics, University of California, Berkeley, CA, USA.}
\email{\href{mailto:gaetz@berkeley.edu}{gaetz@berkeley.edu}.}

\author[Gao]{Yibo Gao}
\address{Beijing International Center for Mathematical Research, Peking University, Beijing, China.}
\email{\href{mailto:gaoyibo@bicmr.pku.edu.cn}{{\tt gaoyibo@bicmr.pku.edu.cn}}}

\author[Jiradilok]{Pakawut Jiradilok}
\address{Department of Mathematics, Massachusetts Institute of Technology, Cambridge, MA, USA.}
\email{\href{mailto:pakawut@mit.edu}{pakawut@mit.edu}.}
\thanks{P.J. was supported by Elchanan Mossel's Vannevar Bush Faculty Fellowship ONR-N00014-20-1-2826 and by Elchanan Mossel's Simons Investigator award (622132).}

\author[Nenashev]{Gleb Nenashev}
\address{Department of Mathematics and Computer Science, St. Petersburg State University, St. Petersburg, Russia.}
\email{\href{mailto:glebnen@gmail.com}{glebnen@gmail.com}.}
\thanks{G.N. was supported by the Knut and Alice Wallenberg Foundation, KAW\,2017.0394} 

\author[Postnikov]{Alexander Postnikov}
\address{Department of Mathematics, Massachusetts Institute of Technology, Cambridge, MA, USA.}
\email{\href{mailto:apost@math.mit.edu}{apost@math.mit.edu}.} 

\date{\today}

\begin{abstract}
We study the maximum multiplicity $\Ans(k,n)$ of a simple transposition $s_k=(k \: k+1)$ in a reduced word for the longest permutation $w_0=n \: n-1 \: \cdots \: 2 \: 1$, a problem closely related to much previous work on sorting networks and on the ``$k$-set" problem.  After reinterpreting the problem in terms of monotone weakly separated paths, we show that, for fixed $k$ and sufficiently large $n$, the optimal density is realized by paths which are periodic in a precise sense, so that  
\[
\Ans(k,n)=c_k n + p_k(n)
\]
for a periodic function $p_k$ and constant $c_k$.  In fact we show that $c_k$ is always rational, and compute several bounds and exact values for this quantity with \emph{repeatable patterns}, which we introduce.
\end{abstract}

\maketitle

\section{Introduction}
\label{sec:intro}

Write $s_k=(k\ k{+}1)$ for the adjacent transpositions in the symmetric group $S_n$.  A \emph{reduced word} for a permutation $w \in S_n$ is an expression $w=s_{i_1}\cdots s_{i_{\ell}}$ of minimal length, and in this case $\ell=\ell(w)$ is called the \emph{length} of $w$; we write $\Red(w)$ for the set of reduced words of $w$.  

There is a unique permutation $w_0=n \: n-1 \: \cdots \: 2 \: 1$ of maximum length ${n \choose 2}$, called the \emph{longest permutation}.  Reduced words of $w_0$ have been extensively studied, as maximal chains in the weak Bruhat order \cite{Edelman-Greene}, in total positivity and cluster algebras, and in the context of random sorting networks \cite{Angel}.  It is not hard to see that the minimum multiplicity of $s_k$ in a reduced word for $w_0$ is $\min(k,n-k)$ (see Section~\ref{sec:other-types}), while the average multiplicity can be computed using the Edelman--Greene bijection \cite{Edelman}.  This paper describes our study of the quantity $\Ans(k,n)$, the \emph{maximum} multiplicity of $s_k$ among all reduced words of $w_0$.  This problem is considerably more difficult, as evidenced by its close connection to the well-known ``$k$-set problem". The maximum multiplicity problem for reduced words of general permutations has been studied by Tenner \cite{tenner2020range}, who gave bounds expressed in terms of permutation patterns.

Throughout much of this paper\footnote{An extended abstract of this work appears in the proceedings of FPSAC \cite{fpsac}.} we consider \emph{monotone weakly separated paths} or \emph{generalized wiring diagrams} instead of reduced words themselves.  From this perspective certain periodicity phenomena appear which are obscured when considering reduced words or their associated pseudoline arrangements.

\subsection{Relation to the $k$-set problem}
\label{sec:intro-k-sets}
Given a collection $A$ of $n$ distinct points in $\mathbb{R}^2$, a \emph{$k$-set} is a subset $B \subseteq A$ of size $k$ which can be separated from $A - B$ by a straight line in $\mathbb{R}^2$.  The \emph{$k$-set problem}, studied since work of Lov\'{a}sz \cite{Lovasz} and Erd\H{o}s--Lov\'{a}sz--Simmons--Straus \cite{Erdos} in the 1970s, asks for the maximum number of $k$-sets admitted by any collection $A$.  This problem has since found application in the analysis of some geometric algorithms.

A common approach to this problem proceeds by first applying projective duality to recast the problem in terms of regions of height $k$ in an arrangement of $n$ lines, and then relaxing it by considering arrangements of $n$ \emph{pseudolines} (curves in the plane such that each pair crosses exactly once).  Many of the strongest known results for the $k$-set problem work with this relaxation, and all available data \cite{Abrego} indicates that the answers in fact agree for lines and for pseudolines.  An arrangement of $n$ pseudolines can equivalently be thought of as the \emph{wiring diagram} for a reduced word of $w_0$, and in this context the problem becomes to maximize the total number of $s_k$'s and $s_{n-k}$'s appearing.  We show in Section~\ref{sec:rationality} that there is a well-defined slope $c_k$ defined by $\Ans(k,n) \sim c_k n$ and that this quantity is the same whether we consider the total multiplicity of $s_k$ and $s_{n-k}$ or just that of $s_k$, so that our original problem is very closely linked to the (pseudoline version of) the $k$-set problem.  

\subsection{Relation to weak separation}\label{sec:intro-weak-separation}
Given a reduced word of a permutation $w$, we can associate a \emph{weakly separated collection} to it, and more specifically, a \emph{monotone weakly separated path}. This process can be viewed as first obtaining a \emph{plabic graph} from the reduced word, and then taking certain face labels. Weakly separated collections are fundamental objects in the theory of the totally nonnegative Grassmannian and related cluster algebras (see, e.g. \cite{Postnikov-icm}). In particular, Oh-Postnikov-Speyer \cite{OPS} constructed a bijection between reduced plabic graphs of any positroid $\mathcal{M}$ and certain maximal weakly separated collections, establishing the purity property. Moreover, maximal weakly separated collections of ${[n]\choose k}$ correspond to \emph{Pl\"ucker clusters} of the Grassmannian $\mathrm{Gr}(k,n)$, which behave nicely among other clusters. In fact, in this paper, we will very often think of a reduced word via its corresponding monotone weakly separated path. We elaborate on this connection in Section~\ref{sec:prelim}.

\subsection{Outline and main results}

In Section~\ref{sec:prelim} we introduce monotone weakly separated paths and establish an equivalent version of the main problem in these terms.  Section~\ref{sec:bounds} introduces \emph{\Gleb diagrams} and applies these to give bounds and some exact values for $\Ans(k,n)$. \TitleGleb diagrams and their weights give a tool for computing upper bounds on $\Ans(k,n)$, while \emph{repeatable patterns}, also introduced in Section~\ref{sec:bounds}, allow explicit constructions of reduced words for all $n$ at once, and thus for determining lower bounds on $\Ans(k,n)$. This technology allows us to show:

\begin{thm}[See Section~\ref{sec:bounds}]
For $k=1,2,3$, the quantity $c_k$ exists and we have $c_1=1, c_2=\frac{3}{2},$ and $c_3=\frac{11}{6}$. Furthermore, explicit reduced words realizing $\mathcal{M}(k,n)$ for $k=1,2,3$ and $n \in \mathbb{N}$ can be obtained from the repeatable patterns given in Section~\ref{sec:bounds}. 
\end{thm}

In Section~\ref{sec:rationality} we introduce \emph{generalized wiring diagrams}, which, for arbitrary fixed $k$, can be used to reason about $\mathcal{M}(k,n)$ for all $n$ simultaneously. We use these objects to show that for all $k$ the quantity
\[
c_k \coloneqq \lim_{n \to \infty} \frac{\Ans(k,n)}{n}
\]
exists, is rational, and is equal to the corresponding limit which counts multiplicities of both $s_k$ and $s_{n-k}$. In fact, what we prove is much stronger:

\begin{thm}[See Section~\ref{sec:rationality}]
For fixed $k$ and sufficiently large $n$, $c_k$ is realized by diagrams which are are periodic in a precise sense, so that computing $c_k$ reduces to a finite search for repeatable patterns. 
\end{thm}  

Finally, in Section~\ref{sec:other-types} we discuss the problem (which is easy for the symmetric group) of \emph{minimizing} the multiplicity of $s_k$ in a reduced word for the longest element $w_0$ in other finite Coxeter groups.

\section{Reduced words and weakly separated paths}\label{sec:prelim}
In this section, we establish relations between reduced words and monotone weakly separated paths. We say that two different sets $I,J\subset[n]$ of cardinality $k$ are \emph{weakly separated} if $\max I- J<\min J- I$ or $\max J- I<\min I- J$, and that a collection of cardinality $k$ subsets of $[n]$ is \emph{weakly separated} if each pair of sets is weakly separated. Note that being weakly separated is not a transitive relation. A sequence of subsets $(A_0,A_1,\ldots,A_N)$ is a \emph{monotone weakly separated path} if the collection $\{A_0,\ldots,A_N\}$ is weakly separated and for each $i=1,\ldots,N$, both $A_{i}- A_{i-1}=:\{x_i\}$ and $A_{i-1}- A_i=:\{y_i\}$ are singleton sets with $x_i>y_i$. 

Given a reduced word $s_{i_1}\cdots s_{i_{\ell}}=\mathbf{i}\in\Red(w)$, and a fixed simple generator $s_k=(k\ k+1)$, let $a_1<\cdots <a_N$ be the positions of all $s_k$'s in $\mathbf{i}$. We obtain permutations $w^{(j)}=s_{i_1}s_{i_2}\cdots s_{i_{a_j}}$ as the products of prefixes of $\mathbf{i}$, where $w^{(0)}=\mathrm{id}$. For $j=1,\ldots,N$, let $A_j=\{w^{(j)}(1),w^{(j)}(2),\ldots,w^{(j)}(k)\}$ be the set of values of $w^{(j)}$ on inputs $1,\ldots,k$, and write $P_k(\mathbf{i})=(A_0,A_1,\ldots,A_N)$. 
\begin{defn}
Given a reduced word $\mathbf{i}=s_{i_1}\cdots s_{i_{\ell}}$ of $w\in S_n$, its corresponding \textit{wiring diagram} consists of wires labeled by $1,2,\ldots,n$ starting at \text{levels} $1,2,\ldots$ respectively from top to bottom, traveling from left to right such that at each timestamp $t$, the two wires at levels $i_t$ and $i_{t+1}$ cross. 
\end{defn}
We will be mainly using wiring diagrams as visualizations for reduced words.
\begin{example}\label{ex:wiring-diagram}
Consider the following reduced word of the longest permutation $w_0\in S_6$:
\[\mathbf{i}=s_3s_2s_1s_3s_2s_3s_4s_3s_5s_4s_3s_2s_1s_3s_2\]
with its corresponding wiring diagram shown in Figure~\ref{fig:wd-123-456}.
\begin{figure}[h!]
\centering
\begin{tikzpicture}[scale=0.6]
\def\a{0.15};
\node[left] at (-1,0) {$1$};
\node[left] at (-1,-1) {$2$};
\node[left] at (-1,-2) {$3$};
\node[left] at (-1,-3) {$4$};
\node[left] at (-1,-4) {$5$};
\node[left] at (-1,-5) {$6$};
\draw(-1,0)--(3-\a,-0)--(3+\a,-1)--(5-\a,-1)--(5+\a,-2)--(6-\a,-2)--(6+\a,-3)--(7-\a,-3)--(7+\a,-4)--(9-\a,-4)--(9+\a,-5)--(16,-5);
\draw(-1,-1)--(2-\a,-1)--(2+\a,-2)--(4-\a,-2)--(4+\a,-3)--(6-\a,-3)--(6+\a,-2)--(8-\a,-2)--(8+\a,-3)--(10-\a,-3)--(10+\a,-4)--(16,-4);
\draw(-1,-2)--(1-\a,-2)--(1+\a,-3)--(4-\a,-3)--(4+\a,-2)--(5-\a,-2)--(5+\a,-1)--(12-\a,-1)--(12+\a,-2)--(14-\a,-2)--(14+\a,-3)--(16,-3);
\draw(-1,-3)--(1-\a,-3)--(1+\a,-2)--(2-\a,-2)--(2+\a,-1)--(3-\a,-1)--(3+\a,0)--(13-\a,0)--(13+\a,-1)--(15-\a,-1)--(15+\a,-2)--(16,-2);
\draw(-1,-4)--(7-\a,-4)--(7+\a,-3)--(8-\a,-3)--(8+\a,-2)--(11-\a,-2)--(11+\a,-3)--(14-\a,-3)--(14+\a,-2)--(15-\a,-2)--(15+\a,-1)--(16,-1);
\draw(-1,-5)--(9-\a,-5)--(9+\a,-4)--(10-\a,-4)--(10+\a,-3)--(11-\a,-3)--(11+\a,-2)--(12-\a,-2)--(12+\a,-1)--(13-\a,-1)--(13+\a,0)--(16,0);
\node at (0,-2.5) {$123$};
\node at (2.5,-2.5) {$124$};
\node at (5,-2.5) {$134$};
\node at (7,-2.5) {$234$};
\node at (9.5,-2.5) {$345$};
\node at (12.5,-2.5) {$346$};
\node at (15,-2.5) {$456$};
\end{tikzpicture}
\caption{The wiring diagram of the reduced word $\mathbf{i}=s_3s_2s_1s_3s_2s_3s_4s_3s_5s_4s_3s_2s_1s_3s_2$}\label{fig:wd-123-456}
\end{figure}
Now fix $k=3$ where $s_k$ appears $6$ times in $\mathbf{i}$. We have the intermediate permutations $w^{(0)}=123456$, $w^{(1)}=124356$, $w^{(2)}=413256$, $w^{(3)}=432156$, $w^{(4)}=435216$, $w^{(5)}=436521$, $w^{(6)}=645321$. Taking their first $k$ values, we obtained $A_0=\{1,2,3\}$, $A_1=\{1,2,4\}$, $A_2=\{1,3,4\}$, $A_3=\{2,3,4\}$, $A_4=\{3,4,5\}$, $A_5=\{3,4,6\}$, $A_6=\{4,5,6\}$ as shown in Figure~\ref{fig:wd-123-456}.
\end{example}
\begin{prop}\label{prop:equivalence-between-words-monotone-path}
Let $P_k(\mathbf{i})$ be constructed as above. Then $P_k(\mathbf{i})$ is a monotone weakly separated path. Conversely, for any monotone weakly separated path $P$ that starts with $\{1,2,\ldots,k\}$, there exists a reduced word $\mathbf{i}$ such that $P_k(\mathbf{i})=P$.
\end{prop}
\begin{proof}
Let $\mathbf{i}\in\Red(w)$ and $P_k(\mathbf{i})=(A_0,\ldots,A_N)$. If some $A_j$ and $A_{j'}$ with $j<j'$ are not weakly separated, then there exists $a\in A_j- A_{j'}$ and $a'\in A_{j'}- A_j$ such that $a>a'$. By definition, $w^{(j)}<w^{(j')}$ in the right weak Bruhat order, but $(a,a')$ is a left inversion of $w^{(j)}$, not of $w^{(j')}$, contradiction. In other words, if we consider the wiring diagram associated to $\mathbf{i}$, the wires labeled $a$ and $a'$ must intersect from $A_0$ to $A_j$, and intersect again from $A_j$ to $A_{j'}$, meaning that $\mathbf{i}$ cannot be reduced. As a result, $\{A_0,\ldots,A_N\}$ is a weakly separated collection. At the same time, $A_j=A_{j-1}-\{x\}\cup\{y\}$ if we write $(x\ y)s_{i_1}\cdots s_{i_{a_j-1}}=s_{i_1}\cdots s_{i_{a_j-1}}s_{i_{a_j}}$. And $x<y$ since $\mathbf{i}$ is reduced. Thus, $P_k(\mathbf{i})=(A_0,\ldots,A_N)$ is a monotone weakly separated path.

Now suppose that we are given a monotone weakly separated path $P=(A_0,\ldots,A_N)$ with $A_0=\{1,\ldots,k\}$. Start with $w^{(0)}=\mathrm{id}$. We are going to construct $w^{(1)},w^{(2)},\ldots$ with a reduced word $\mathbf{i}$ along the way such that $P_k(\mathbf{i})=P$. Suppose that we have constructed $w^{(j)}=s_{i_1}\cdots s_{i_{m}}$ and let $x\in A_j- A_{j+1}$, $y\in A_{j+1}- A_j$ with $x<y$. Suppose that $w^{(j)}(a)=x$ and $w^{(j)}(b)=y$ with $a\leq k<b$. We can continue the construction of $\mathbf{i}$ by $w^{(j+1)}=w^{(j)} (s_{a}s_{a+1}\cdots s_{k-1})(s_{b-1}s_{b-2}\cdots s_{k+1})s_k$. Here, $s_{a}s_{a+1}\cdots s_{k-1}$ moves $x$ from position $a$ to position $k$ while $s_{b-1}s_{b-2}\cdots s_{k+1}$ moves $y$ from position $b$ to position $k+1$. In the end, the $s_k$ exchanges the values $x$ and $y$. Therefore, we automatically have $\{w^{(j+1)}(1),\ldots,w^{(j+1)}(k)\}=A_j-\{x\}\cup\{y\}=A_{j+1}$ as desired. The only thing left to show is that the word $\mathbf{i}$ coming from such construction is reduced.

If $\mathbf{i}$ is not reduced, we can without loss of generality assume that in some step when we are constructing $w^{(j+1)}$ from $w^{(j)}$, a simple generator $s_p$ exchanges a larger value at position $p$ with a smaller value at position $p+1$. Keep the notation as in the above paragraph. We cannot have $p=k$ since $s_k$ always exchanges $A_{j}- A_{j+1}$ at position $k$ with $A_{j+1}- A_j$ at position $k+1$. So by symmetry, we assume $p<k$, and that such $s_p$ exchanges value $x\in A_{j+1}- A_j$ at position $p$ with value $z$ at position $p+1$, with $x>z$. Since $z<x$, the values $z$ and $x$ must have been switched before, when we are constructing $w^{(j'+1)}$ from $w^{(j')}$, with $j'<j$. By construction, we are either moving $z$ out of $A_{j'}$ to $A_{j'+1}$, or moving $x$ into $A_{j'+1}$ from out of $A_{j'}$. In both cases, $z\notin A_{j'+1}$ and $x\in A_{j'+1}$. As a result, $x\in A_{j'+1}- A_{j+1}$, $z\in A_{j+1}- A_{j'+1}$, but $z<x$. As $A_{j+1}$ and $A_{j'+1}$ are weakly separated, we must have $\max A_{j+1}- A_{j'+1}<\min A_{j'+1}- A_{j+1}$. But $j'<j$, there cannot possibly be a monotone path from $A_{j'+1}$ to $A_{j+1}$. Contradiction. Thus, this construction results in a reduced word $\mathbf{i}$ as desired.
\end{proof}

Consequently, we say that $P_k(\mathbf{i})$ is the monotone weakly separated path associated to $\mathbf{i}\in\Red(w)$. Clearly, if $P_k(\mathbf{i})$ consists of $N+1$ subsets from $A_0$ to $A_N$, then there are exactly $N$ $s_k$'s in $\mathbf{i}$. Proposition~\ref{prop:equivalence-between-words-monotone-path} allows us to translate the problem of finding the maximal number of $s_k$'s in $\Red(w)$ to finding the longest monotone weakly separated path that starts at $\{1,2,\ldots,k\}$.

\section{Repeatable patterns and \Gleb diagrams}\label{sec:bounds}
This section introduces \Gleb diagrams and repeatable patterns, and shows:

\begin{itemize}
\item[(i)] $\mathcal{M}(1,n) = n-1$, for every integer $n \ge 2$, 
\item[(ii)] (Theorem~\ref{thm:M2n}) $\mathcal{M}(2,n) = \left\lceil \frac{3}{2} n \right\rceil - 3$, for every integer $n \ge 3$, and 
\item[(iii)] (Theorem~\ref{thm:M3n}) $\mathcal{M}(3,n) = \left\lceil \frac{11}{6} n \right\rceil - 5$, for every integer $n \ge 4$.
\end{itemize}

\subsection{\TitleGleb diagrams}\label{ss:Gd} 
Let $k$ and $n$ be positive integers such that $1 \le k \le n-1$. Suppose that $P = (A_0, A_1, \ldots, A_N)$ is a monotone weakly separated path from $A_0 = \{1, 2, \ldots, k\}$ to $A_N = \{n-k+1, \ldots, n\}$. We define the {\em \Gleb diagram} $\mathbf{D}(P)$ of $P = (A_0, A_1, \ldots, A_N)$ to be the simple undirected graph on the vertex set $[n] = \{1, 2, \ldots, n\}$ in which an edge $(i,j)$ appears if and only if there exists $a \in [N]$ such that $\{i, j\} = (A_{a-1} - A_a) \cup (A_a - A_{a-1})$.

We give a quick remark about the above definition of $\mathbf{D}(P)$. Note that, for each pair $\{i,j\}$, if there exists $a \in [N]$ such that $\{i,j\} = (A_{a-1} - A_a) \cup (A_a - A_{a-1})$, then such an index $a$ is unique: the equation implies that in the corresponding wiring diagram, the $a^{\text{th}}$ $s_k$-crossing from the left is a crossing between wire $i$ and wire $j$, and any two wires intersect exactly once in a wiring diagram of a reduced word of $w_0$. (See Figures~\ref{fig:wd-123-456} and \ref{fig:Gleb-123-456} for an example.)

\begin{example}\label{ex:G-wd}
When $n = 6$ and $k = 3$, an example of a monotone weakly separated path is $P = 123 - 124 - 134 - 234 - 345 - 346 - 456$, realized in Example~\ref{ex:wiring-diagram}. Here, the shorthand $ijk$ represents the triple $\{i,j,k\}$. The \Gleb diagram $\mathbf{D}(P)$ of $P$ is shown in Figure~\ref{fig:Gleb-123-456}.
\end{example}
\begin{figure}[h!]
\begin{center}
\def\ptr{2}
\def\h{0.5}
\def\glebarc[#1,#2]{
\draw[line width = 1 pt] (#1,0) arc(180:0:{0.5*#2} and {0.5*#2});}
    \begin{tikzpicture}
    
    \filldraw[black] (1,0) circle (\ptr pt);
    \filldraw[black] (2,0) circle (\ptr pt);
    \filldraw[black] (3,0) circle (\ptr pt);
    \filldraw[black] (4,0) circle (\ptr pt);
    \filldraw[black] (5,0) circle (\ptr pt);
    \filldraw[black] (6,0) circle (\ptr pt);
    
    \node at (1,0-\h) {$1$};
    \node at (2,0-\h) {$2$};
    \node at (3,0-\h) {$3$};
    \node at (4,0-\h) {$4$};
    \node at (5,0-\h) {$5$};
    \node at (6,0-\h) {$6$};
    
    \glebarc[1,1]
    \glebarc[2,1]
    \glebarc[3,1]
    \glebarc[5,1]
    \glebarc[3,2]
    \glebarc[2,3]
    
    \end{tikzpicture}
\end{center}
\caption{The \Gleb diagram $\mathbf{D}(P)$ of the monotone weakly separated path $P$ in Example~\ref{ex:G-wd}.}\label{fig:Gleb-123-456}
\end{figure}

Given an \Gleb diagram, we put each vertex $i \in [n]$ of the diagram at the point $(i,0) \in \mathbb{R}^2$ and draw each edge $(i,j)$ as a semicircle on the upper-half plane. Imagine that each semicircular curve in the \Gleb diagram has weight $1$. For each curve, assume that the weight is distributed uniformly across the horizontal length ({\em not} the curve length). For example, if an edge $e$ joins $(1,0)$ and $(4,0)$, then there is weight $2/3$ {\em above} $[2,4]$ coming from $e$. If we have a finite collection of curves, define the total weight as the sum of individual weights. Note that $\mathcal{M}(k,n)$ is the maximum possible total weight in an \Gleb diagram.

\begin{proposition}\label{prop:upper-bound-series}
For any positive integers $k$ and $n$ with $1 \le k \le n-1$, we have
\begin{equation}\label{eq:upper-bound-series}
\mathcal{M}(k,n) \le \Bigg( \underbrace{1 + \frac{1}{2} + \frac{1}{2} + \frac{1}{3} + \frac{1}{3} + \frac{1}{3} + \frac{1}{4} + \frac{1}{4} + \frac{1}{4} + \frac{1}{4} + \cdots }_{k \text{ terms}}\Bigg) \cdot n.
\end{equation}
\end{proposition}

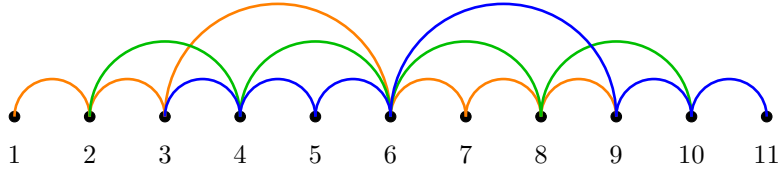
\begin{figure}
\begin{center}
\def\ptr{2}
\def\h{0.5}
\def\glebarccol[#1,#2,#3]{
\draw[color = #3, line width = 1 pt] (#1,0) arc(180:0:{0.5*#2} and {0.5*#2});}
    \begin{tikzpicture}
    
    \filldraw[black] (1,0) circle (\ptr pt);
    \filldraw[black] (2,0) circle (\ptr pt);
    \filldraw[black] (3,0) circle (\ptr pt);
    \filldraw[black] (4,0) circle (\ptr pt);
    \filldraw[black] (5,0) circle (\ptr pt);
    \filldraw[black] (6,0) circle (\ptr pt);
    \filldraw[black] (7,0) circle (\ptr pt);
    \filldraw[black] (8,0) circle (\ptr pt);
    \filldraw[black] (9,0) circle (\ptr pt);
    \filldraw[black] (10,0) circle (\ptr pt);
    \filldraw[black] (11,0) circle (\ptr pt);
    
    \node at (1,0-\h) {$1$};
    \node at (2,0-\h) {$2$};
    \node at (3,0-\h) {$3$};
    \node at (4,0-\h) {$4$};
    \node at (5,0-\h) {$5$};
    \node at (6,0-\h) {$6$};
    \node at (7,0-\h) {$7$};
    \node at (8,0-\h) {$8$};
    \node at (9,0-\h) {$9$};
    \node at (10,0-\h) {$10$};
    \node at (11,0-\h) {$11$};
    
    \glebarccol[1,1,orange]
    \glebarccol[2,1,orange]
    \glebarccol[3,3,orange]
    \glebarccol[6,1,orange]
    \glebarccol[7,1,orange]
    \glebarccol[8,1,orange]
    \glebarccol[2,2,green!75!black]
    \glebarccol[4,2,green!75!black]
    \glebarccol[6,2,green!75!black]
    \glebarccol[8,2,green!75!black]
    \glebarccol[3,1,blue]
    \glebarccol[4,1,blue]
    \glebarccol[5,1,blue]
    \glebarccol[6,3,blue]
    \glebarccol[9,1,blue]
    \glebarccol[10,1,blue]
    
    \end{tikzpicture}
\end{center}
\caption{An example of the decoration performed in the proof of Proposition~\ref{prop:upper-bound-series}. The diagram above is the decorated \Gleb diagram of the monotone weakly separated path $P = \big( \{1,2,3\},\{1,2,4\},\{1,2,5\},\{1,4,5\},\{2,4,5\},\{3,4,5\},\{4,5,6\}$, $\{4,5,7\},\{5,6,7\},\{5,7,8\},\{6,7,8\},\{7,8,9\},\{7,8,10\},\{7,8,11\}$, $\{7,10,11\},\{8,10,11\},\{9,10,11\}\big)$.}\label{fig:dec-Gleb}
\end{figure}

\begin{proof}
Let $P$ be a monotone weakly separated path from $A_0 = \{1, 2, \ldots, k\}$ to $A_N = \{n-k+1, \ldots, n\}$. Perform the following auxiliary decoration using $k$ different colors $\col_1, \col_2, \ldots, \col_k$. First, we color $1, 2, \ldots, k$ in $A_0$ so that $i$ gets color $\col_i$. Each time we go from $A_j$ to $A_{j+1}$, if we have $A_{j+1} = A_j \cup \{y\} - \{x\}$, then color $y$ in $A_{j+1}$ the same color as $x$ in $A_j$. Also color the semicircle connecting $(x,0)$ and $(y,0)$ with the same color that we used to color $x \in A_j$ and $y \in A_{j+1}$. For each number $z \in A_j \cap A_{j+1}$, color $z$ in $A_{j+1}$ the same color as we color $z$ in $A_j$. See Figure~\ref{fig:dec-Gleb} for an example.

As a result of the decoration, the $k$ numbers in $A_N$ are now colored with $k$ different colors. There is a permutation $\sigma \in S_k$ such that for each $i \in [k]$, the number $n-k+\sigma(i)$ has color $\col_i$. Now there are $k$ continuous curves $\gamma_1, \gamma_2, \ldots, \gamma_k$ such that $\gamma_i$ has endpoints $(i,0)$ and $(n-k+\sigma(i),0)$, and has color $\col_i$ from one end to the other.

From the coloring argument above, there are at most $k$ pieces of different semicircles in the vertical strip above $[i,i+1]$ for each $i \in [n-1]$. Furthermore, for each $t \in \mathbb{Z}_{\ge 1}$, there exist at most $t$ pieces that come from semicircles whose diameters are exactly $t$. Each such piece contributes the weight of $1/t$. Therefore, the weight of the \Gleb diagram above $[i,i+1]$ is at most $1 + 1/2 + 1/2 + 1/3 + 1/3 + 1/3 + \cdots$, and hence the total weight is at most the right-hand side of (\ref{eq:upper-bound-series}). 
\end{proof}

By estimating the summation in Proposition~\ref{prop:upper-bound-series}, we obtain the following corollary.

\begin{corollary}\label{cor:sqrt(2k)}
For any positive integers $k$ and $n$ such that $1 \le k \le n-1$, we have the upper bound $\mathcal{M}(k,n) \le \sqrt{2k} \cdot n$.
\end{corollary}

\begin{remark}
Together with Theorem~\ref{thm:dk} below, our arguments above give a short proof of the upper bound $O(\sqrt{k}) \cdot n$ for the pseudoline $k$-sets problem. It might be instructive to compare this bound with the known upper bounds in the literature of the straight line setting. For the classical planar $k$-sets problem, Pach, Steiger, and Szemer\'{e}di \cite{PSS92} have shown the upper bound
\begin{equation}\label{PSS92-log-star}
O\!\left( \frac{\sqrt{k}}{\log^* k} \right) \cdot n,
\end{equation}
which is slightly stronger than $O(\sqrt{k}) \cdot n$. Our pseudoline setting is more general. We do not know if the upper bound~\eqref{PSS92-log-star} of Pach--Steiger--Szemer\'{e}di holds for $\mathcal{M}(k,n)$ or not. This might be an interesting direction to further investigate.
\end{remark}

The upper bound in Proposition~\ref{prop:upper-bound-series} can be slightly improved as follows. First, note that the number of unit segments $[i,i+1]$ is actually $n-1$ (instead of $n$). Second, note that the segments $[i,i+1]$ near the {\em ends} (vertices $(1,0)$ and $(n,0)$) should have smaller upper bounds because there are fewer than $k$ pieces of curve above those segments. This improvement leads to sharp results when $k$ is small. For $k = 1$, we find that $\mathcal{M}(1,n) \le n-1$, for every $n \ge 2$. For $k = 2$, we find $\mathcal{M}(2,n) \le \frac{3n-5}{2}$, for every $n \ge 3$. Since $\mathcal{M}(2,n)$ is an integer, we can write $\mathcal{M}(2,n) \le \left\lceil \frac{3}{2} n \right\rceil - 3$. We will see in the next subsection that these bounds in the cases $k = 1$ and $k = 2$ are sharp.

\subsection{$k \le 2$ and repeatable patterns}\label{ss:RP}

In the previous subsection, we have seen that $\mathcal{M}(1,n) \le n-1$ for every $n \ge 2$. In fact, it is easy to see that $\mathcal{M}(1,n) = n-1$. Indeed, the sequence $(\{1\}, \{2\}, \ldots, \{n\})$ is a monotone weakly separated path.

Things get more interesting when $k = 2$. In the previous subsection, we have also seen that $\mathcal{M}(2,n) \le \left\lceil \frac{3}{2} n \right\rceil - 3$, for each $n \ge 3$. Now we claim that the inequality is in fact an equality by giving explicit constructions using the idea of {\em repeatable patterns}.

In the definition below, if $S$ is a finite set of integers and $t$ is an integer, we write $S+t$ to denote $\{s + t : s \in S\}$.

\begin{defn}\label{defn:rp}
Let $L$ and $d$ be positive integers. A {\em repeatable pattern} $R$ with {\em parameters} $(L, d)$ is a monotone weakly separated path $R = (A_0, A_1, \ldots, A_L)$ which satisfies the following conditions:
\begin{itemize}
\item $A_L = A_0 + d$, and
\item for any positive integer $m$, the sequence
\begin{align*}
\Big( A_0, &A_1, A_2, \ldots, A_L, \\
&A_1 + d, A_2 + d, \ldots, A_L + d, \\
&A_1 + 2d, A_2 + 2d, \ldots, A_L + 2d, \\
& \vdots \\
&A_1 + md, A_2 + md, \ldots, A_L + md \Big)
\end{align*}
is a monotone weakly separated path.
\end{itemize}
\end{defn}

\begin{example}
The pattern $12 - 13 - 23 - 34$ is a repeatable pattern with parameters $(L, d) = (3, 2)$. Here, we use the shorthand $ab$ to denote $\{a,b\}$. By concatenation, the pattern gives the infinite sequence
\[
12 - 13 - 23 - 34 - 35 - 45 - 56 - 57 - 67 - 78 - 79 - 89 - \cdots.
\]
Any finite prefix of the infinite sequence above is a monotone weakly separated path.
\end{example}

\begin{thm}\label{thm:M2n}
For each positive integer $n \ge 3$, we have $\mathcal{M}(2,n) = \left\lceil \frac{3}{2} n \right\rceil - 3$.
\end{thm}
\begin{proof}
The first $\left\lceil \frac{3}{2} n \right\rceil - 2$ terms of the infinite sequence in the previous example is a monotone weakly separated path from $\{1, 2\}$ to $\{n-1, n\}$. Combine this construction with the upper bound for $\mathcal{M}(2,n)$ above to finish.
\end{proof}

For general $k$, the existence of a repeatable pattern yields a lower bound for $\mathcal{M}(k,n)$. In Section~\ref{sec:rationality}, we will see that for every positive integer $k$, the limit $c_k := \lim_{n \to \infty} \frac{\mathcal{M}(k,n)}{n}$ exists. The existence of a repeatable pattern $R$ with parameters $(L, d)$ immediately yields the lower bound $c_k \ge \frac{L}{d}$. It turns out, as Theorem~\ref{thm:ck-is-rational} below shows, that $c_k = \max_R L/d$, where the maximization is over all repeatable patterns $R = (A_0, A_1, \ldots, A_L)$ with $|A_0| = k$, and the ratio $L/d$ depends on the repeatable pattern $R$. In particular, the maximum (not just the supremum) exists. 

\subsection{\TitleGleb diagrams when $k = 3$, part i: decomposition algorithm}\label{ss:D}

We use the ``decomposition algorithm" (Algorithm~\ref{alg:fish}) below to break the interval $[1,n]$ in the \Gleb diagram into smaller non-overlapping intervals in a way that we can prove upper bounds of weights for these intervals separately.

We will show that $\mathcal{M}(3,n) = \left\lceil \frac{11}{6} n \right\rceil - 5$, for $n \ge 4$. The cases $n=4$ and $n=5$ can be readily taken care of. By using $\mathcal{M}(k,n) = \mathcal{M}(n-k,n)$, we find that $\mathcal{M}(3,4) = \mathcal{M}(1,4) = 3$ and $\mathcal{M}(3,5) = \mathcal{M}(2,5) = 5$. For the rest of this subsection, assume $n \ge 6$.

Suppose that an \Gleb diagram coming from a monotone weakly separated path from $\{1, 2, 3\}$ to $\{n-2, n-1, n\}$ is given. Write $\wt$ to denote the weight function, so that if $I \subseteq [1,n]$ is an interval, then $\wt(I)$ is the weight above $I$. 

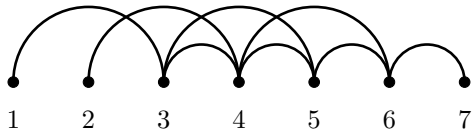
\begin{figure}
\begin{center}
\def\ptr{2}
\def\h{0.5}
\def\glebarc[#1,#2]{
\draw[line width = 1 pt] (#1,0) arc(180:0:{0.5*#2} and {0.5*#2});}
    \begin{tikzpicture}
    
    \filldraw[black] (1,0) circle (\ptr pt);
    \filldraw[black] (2,0) circle (\ptr pt);
    \filldraw[black] (3,0) circle (\ptr pt);
    \filldraw[black] (4,0) circle (\ptr pt);
    \filldraw[black] (5,0) circle (\ptr pt);
    \filldraw[black] (6,0) circle (\ptr pt);
    \filldraw[black] (7,0) circle (\ptr pt);
    
    \node at (1,0-\h) {$1$};
    \node at (2,0-\h) {$2$};
    \node at (3,0-\h) {$3$};
    \node at (4,0-\h) {$4$};
    \node at (5,0-\h) {$5$};
    \node at (6,0-\h) {$6$};
    \node at (7,0-\h) {$7$};
    
    \glebarc[3,1]
    \glebarc[4,1]
    \glebarc[5,1]
    \glebarc[6,1]
    \glebarc[1,2]
    \glebarc[2,2]
    \glebarc[3,2]
    \glebarc[4,2]
    
    \end{tikzpicture}
\end{center}
\caption{The \Gleb diagram $\mathbf{D}(P)$ of the monotone weakly separated path $P$ in Example~\ref{ex:wt-2}.}\label{fig:wt-2-Gleb}
\end{figure}

\begin{example}\label{ex:wt-2}
Consider the following monotone weakly separated path
\[
P = 123 - 124 - 234 - 245 - 246 - 247 - 267 - 467 - 567.
\]
The \Gleb diagram of $P$ is shown in Figure~\ref{fig:wt-2-Gleb}. The weights of the unit intervals of this \Gleb diagram are as follows: $\wt([1,2]) = 1/2$, $\wt([2,3]) = 1$, $\wt([3,4]) = 2$, $\wt([4,5]) = 2$, $\wt([5,6]) = 3/2$, and $\wt([6,7]) = 1$.
\end{example}

We also define the {\em weight limit} function $\wtlim$ as follows. Declare $\wtlim([1,2]) = 1$, $\wtlim([2,3]) = 3/2$, $\wtlim([n-2,n-1]) = 3/2$, and $\wtlim([n-1,n]) = 1$. If $3 \le i \le n-3$, we declare $\wtlim([i,i+1]) = 11/6$. The weight limit function is also defined to satisfy the usual additivity condition: $\wtlim(A \cup B) = \wtlim(A) + \wtlim(B)$ if $A \cap B$ contains no nontrivial interval.

Observe that for a unit interval $[i,i+1] \subseteq [1,n]$ (with $i \in \mathbb{Z}$) to exceed its weight limit, the only possible way is to have $\wt([i,i+1]) = 2$. Moreover, the unit intervals $[1,2]$, $[2,3]$, $[n-2,n-1]$, and $[n-1,n]$ never exceed their weight limits. These observations follow from the definition.

The ``decomposition algorithm" (Algorithm~\ref{alg:fish}) is given below. The input of the algorithm is an \Gleb diagram $\mathbf{D}$ that comes from a monotone weakly separated path from $\{1, 2, 3\}$ to $\{n-2, n-1, n\}$. The output of the algorithm is the collection $\mathcal{F} = \mathcal{F}(\mathbf{D})$ of intervals with non-overlapping interiors.

\begin{algorithm}[decomposition algorithm]\label{alg:fish} \hphantom{=}

\medskip

\noindent \textbf{Input:} an \Gleb diagram $\mathbf{D}$ that comes from a monotone weakly separated path from $\{1, 2 , 3\}$ to $\{n-2, n-1, n\}$.

\medskip

\noindent Begin with an empty collection $\mathcal{F} = \varnothing$.

\medskip

\noindent If for every $i \in [n-1]$, the unit interval $[i,i+1]$ satisfies $\wt([i,i+1]) \le 11/6$, 

\setlength{\leftskip}{25 pt} \noindent then \textbf{output} $\mathcal{F}(\mathbf{D}) = \varnothing$, and we finish the algorithm.

\medskip

\setlength{\leftskip}{0 pt} \noindent On the other hand, if some unit interval exceeds its weight limit, 

\setlength{\leftskip}{25 pt} \noindent then write
\[
[1,n] = L_0 \cup H_1 \cup L_1 \cup H_2 \cup L_2 \cup \cdots \cup L_m,
\]
where $L_0 = [1, \ell_0]$, $H_1 = [\ell_0,h_1]$, $L_1 = [h_1, \ell_1]$, $\ldots$, $L_m = [h_m,n]$, where \\ $1, \ell_0, h_1, \ell_1, \ldots, h_m, n$ is a strictly increasing sequence of positive integers such that every unit interval in any $L_i$ has weight under or equal to its weight limit, and every unit interval in any $H_i$ has weight exceeding its weight limit.

\medskip

\noindent For $i = 1, 2, \ldots, m$:

\medskip

\setlength{\leftskip}{50 pt} \noindent consider the interval $H_i$. We know from Proposition~\ref{prop:mu(H_i)-le-4} (proved below) that $\mu(H_i) \in \{1, 2, 3, 4\}$.

\medskip

\noindent \underline{Case 1.} $\mu(H_i) = 4$. Write $H_i = [a,a+4]$. Add the intervals $[a-1,a+2]$ and $[a+2, a+5]$ into the collection $\mathcal{F}$.

\medskip

\noindent \underline{Case 2.} $\mu(H_i) = 3$. Write $H_i = [a,a+3]$. 

\medskip

\setlength{\leftskip}{75 pt} \noindent \underline{Case 2.1.} There is no semicircle connecting $(a-1,0)$ and $(a,0)$ in $\mathbf{D}$. Add the interval $[a-1,a+3]$ to $\mathcal{F}$.

\medskip

\noindent \underline{Case 2.2.} There is a semicircle connecting $(a-1,0)$ and $(a,0)$ in $\mathbf{D}$. Add the interval $[a,a+4]$ to $\mathcal{F}$.

\medskip

\setlength{\leftskip}{50 pt} \noindent \underline{Case 3.} $\mu(H_i) = 2$. Write $H_i = [a,a+2]$.

\medskip

\setlength{\leftskip}{75 pt} \noindent \underline{Case 3.1.} There is no semicircle connecting $(a-1,0)$ and $(a,0)$ in $\mathbf{D}$. Add the interval $[a-1,a+2]$ to $\mathcal{F}$.

\medskip

\noindent \underline{Case 3.2.} There is a semicircle connecting $(a-1,0)$ and $(a,0)$, but there is no semicircle connecting $(a+2,0)$ and $(a+3,0)$ in $\mathbf{D}$. Add the interval $[a,a+3]$ to $\mathcal{F}$.

\medskip

\noindent \underline{Case 3.3.} There is a semicircle connecting $(a-1,0)$ and $(a,0)$, and also there is a semicircle connecting $(a+2,0)$ and $(a+3,0)$ in $\mathbf{D}$. Add the intervals $[a-2,a+1]$ and $[a+1, a+4]$ to $\mathcal{F}$.

\medskip

\setlength{\leftskip}{50 pt} \noindent \underline{Case 4.} $\mu(H_i) = 1$. Write $H_i = [a,a+1]$.

\medskip

\setlength{\leftskip}{75 pt} \noindent \underline{Case 4.1.} There is no semicircle connecting $(a-1,0)$ and $(a,0)$ in $\mathbf{D}$. Add the interval $[a-1,a+1]$ to $\mathcal{F}$.

\medskip

\noindent \underline{Case 4.2.} There is a semicircle connecting $(a-1,0)$ and $(a,0)$, but there is no semicircle connecting $(a+1,0)$ and $(a+2,0)$ in $\mathbf{D}$. Add the interval $[a,a+2]$ to $\mathcal{F}$.

\medskip

\noindent \underline{Case 4.3.} There is a semicircle connecting $(a-1,0)$ and $(a,0)$, and also there is a semicircle connecting $(a+1,0)$ and $(a+2,0)$ in $\mathbf{D}$. If $\wt([a-2,a+1]) \le \wtlim([a-2,a+1])$, then add the interval $[a-2,a+1]$ to $\mathcal{F}$. If $\wt([a-2,a+1]) > \wtlim([a-2,a+1])$, then add the interval $[a,a+3]$ to $\mathcal{F}$.

\medskip

\setlength{\leftskip}{25 pt} \noindent \textbf{Output} $\mathcal{F}(\mathbf{D}) = \mathcal{F}$, and we finish the algorithm.
\end{algorithm}

\begin{thm}[decomposition theorem]\label{thm:alg-good}
Let $n \ge 6$ be a positive integer. Let $\mathbf{D}$ be an \Gleb diagram of a monotone weakly separated path from $\{1, 2, 3\}$ to $\{n-2, n-1, n\}$. Let $\mathcal{F}(\mathbf{D})$ be the collection of intervals obtained from Algorithm~\ref{alg:fish}. Then
\begin{itemize}
\item[(a)] any two different intervals in $\mathcal{F}(\mathbf{D})$ are non-overlapping, and
\item[(b)] any interval $I \in \mathcal{F}(\mathbf{D})$ satisfies $\wt(I) \le \wtlim(I)$.
\end{itemize}
\end{thm}

The proof of Theorem~\ref{thm:alg-good} will be given in Section~\ref{ss:GDC}. 

\begin{corollary}\label{cor:M3n-upper}
For any positive integer $n \ge 6$, we have
\[
\mathcal{M}(3,n) \le \left\lceil \frac{11}{6} n \right\rceil - 5.
\]
\end{corollary}
\begin{proof}
Take any \Gleb diagram $\mathbf{D}$ of a monotone weakly separated path from $\{1, 2, 3\}$ to $\{n-2, n-1, n\}$ with the maximum possible weight so that $\mathcal{M}(3,n) = \wt([1,n])$. From Theorem~\ref{thm:alg-good}, we have
\begin{align*}
\mathcal{M}(3,n) = \wt([1,n]) &= \wt\!\left( [1,n] - \cup\mathcal{F} \right) + \sum_{I \in \mathcal{F}} \wt(I) \\
&\le \wtlim\!\left( [1,n] - \cup\mathcal{F} \right) + \sum_{I \in \mathcal{F}} \wtlim(I)  \\
&= \wtlim([1,n]) = 1 + \frac{3}{2} + (n-5) \cdot \frac{11}{6} + \frac{3}{2} + 1 = \frac{11n-25}{6}.
\end{align*}
Since $\mathcal{M}(3,n) \in \mathbb{Z}$, we have that $\mathcal{M}(3,n) \le \left\lfloor \frac{11n-25}{6} \right\rfloor = \left\lceil \frac{11}{6} n \right\rceil - 5$, as desired.
\end{proof}

\subsection{\TitleGleb diagrams when $k = 3$, part ii: \Gleb diagram chasing}\label{ss:GDC}
Below we define a useful object called the {\em bicolored \Gleb diagram} $\mathbf{BiD}(P)$. By looking at edges in $\mathbf{BiD}(P)$, we are able to rule out some configurations of edges in the original \Gleb diagram $\mathbf{D}(P)$, a process we call {\em \Gleb diagram chasing}. Using \Gleb diagram chasing, we prove Proposition~\ref{prop:mu(H_i)-le-4}, Lemmas~\ref{l:mu(H_i)=1}, \ref{l:mu(H_i)=1-deg}, \ref{l:mu(H_i)=2}, \ref{l:mu(H_i)=2-deg}, and Proposition~\ref{prop:lengths} which are then used in the proof of Theorem~\ref{thm:alg-good}.

\begin{defn}\label{defn:BiGP}
Let $P$ be a monotone weakly separated path $P = (A_0, A_1, \ldots, A_N)$ from $A_0 = \{1, 2, 3\}$ to $A_N = \{n-2, n-1, n\}$. The {\em bicolored \Gleb diagram} $\mathbf{BiD}(P)$ is the multigraph on the vertex set $[n] = \{1, 2, \ldots, n\}$ together with the coloring $\mathfrak{c}: E\!\left( \mathbf{BiD}(P) \right) \to \left\{ \text{black}, \text{red} \hspace{0.1 em} \right\}$ on the edges defined as follows. The black edges are precisely the edges in the (original) \Gleb diagram $\mathbf{D}(P)$. The red edges are added sequentially. For each $i \in [N]$, let $C_i$ denote the pair $(A_i - A_{i-1}) \cup (A_{i-1} - A_i)$. We add the red edges in $N-1$ steps. In the $j^{\text{th}}$ step, consider the two pairs $C_j$ and $C_{j+1}$. Suppose that $C_j = \{a,b\}$ and $C_{j+1} = \{c,d\}$, with $a < b$ and $c < d$. If $b \neq c$, add a red edge joining $b$ and $c$. If $a \neq d$, add a red edge joining $a$ and $d$.
\end{defn}

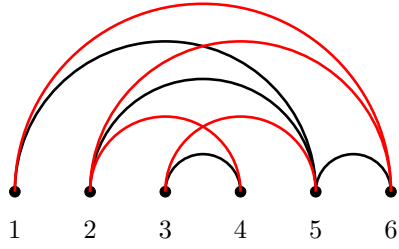
\begin{figure}
\begin{center}
\def\ptr{2}
\def\h{0.5}
\def\glebarccol[#1,#2,#3]{
\draw[color = #3, line width = 1 pt] (#1,0) arc(180:0:{0.5*#2} and {0.5*#2});}
    \begin{tikzpicture}
    
    \filldraw[black] (1,0) circle (\ptr pt);
    \filldraw[black] (2,0) circle (\ptr pt);
    \filldraw[black] (3,0) circle (\ptr pt);
    \filldraw[black] (4,0) circle (\ptr pt);
    \filldraw[black] (5,0) circle (\ptr pt);
    \filldraw[black] (6,0) circle (\ptr pt);
    
    \node at (1,0-\h) {$1$};
    \node at (2,0-\h) {$2$};
    \node at (3,0-\h) {$3$};
    \node at (4,0-\h) {$4$};
    \node at (5,0-\h) {$5$};
    \node at (6,0-\h) {$6$};
    
    \glebarccol[1,4,black]
    \glebarccol[2,3,black]
    \glebarccol[3,1,black]
    \glebarccol[5,1,black]
    \glebarccol[2,2,red]
    \glebarccol[3,2,red]
    \glebarccol[1,5,red]
    \glebarccol[2,4,red]
    
    \end{tikzpicture}
\end{center}
\caption{The bicolored \Gleb diagram $\mathbf{BiD}(P)$ of the monotone weakly separated path $P$ in Example~\ref{ex:BiG-123-456}.}\label{fig:BiG-123-456}
\end{figure}

\begin{example}\label{ex:BiG-123-456}
The bicolored \Gleb diagram $\mathbf{BiD}(P)$ of $P = 123 - 124 - 145 - 146 - 456$ is shown in Figure~\ref{fig:BiG-123-456}.
\end{example}

\begin{proposition}\label{prop:simple-BiGP}
The multigraph $\mathbf{BiD}(P)$ is simple. In other words, each pair of different nodes $i, j \in [n]$ are either (i) joined by one black edge, (ii) joined by one red edge, or (iii) not adjacent.
\end{proposition}

Before proving the proposition, we show a lemma about black edges.

\begin{lem}\label{l:consecutive-curves}
Let $P = (A_0, A_1, \ldots, A_N)$ and $C_1, C_2, \ldots, C_N$ be defined as in Definition~\ref{defn:BiGP}. Suppose that $j$ is a positive integer such that $1 \le j \le N-1$. Write $C_j = \{a, b\}$ and $C_{j+1} = \{c, d\}$, where $a < b$ and $c < d$. Then one of the following six outcomes happens:
\begin{itemize}
\item[(i)] $a < b = c < d$, 
\item[(ii)] $c < d = a < b$,
\item[(iii)] $a < c < b < d$,
\item[(iv)] $c < a < d < b$,
\item[(v)] $a < c < d < b$,
\item[(vi)] $c < a < b < d$.
\end{itemize}
\end{lem}

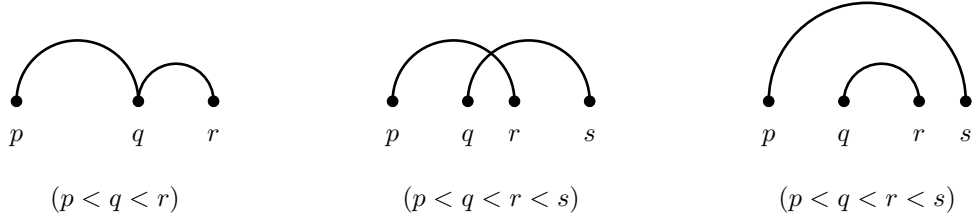
\begin{figure}
\begin{center}
\def\ptr{2}
\def\h{0.5}
\def\glebarc[#1,#2]{
\draw[line width = 1 pt] (#1,0) arc(180:0:{0.5*#2} and {0.5*#2});}
    \begin{tikzpicture}
    \begin{scope}[shift={(-5,0)}]
    \filldraw[black] (0,0) circle (\ptr pt);
    \filldraw[black] (1.62,0) circle (\ptr pt);
    \filldraw[black] (2.62,0) circle (\ptr pt);
    
    \node[anchor = north] at (0,0-0.5*\h) {$p$};
    \node[anchor = north] at (1.62,0-0.5*\h) {$q$};
    \node[anchor = north] at (2.62,0-0.5*\h) {$r$};
    
    \glebarc[0,1.62]
    \glebarc[1.62,1]
    
    \node[anchor = north] at (1.31,0-2*\h) {$(p<q<r)$};
    \end{scope}
    
    \begin{scope}[shift={(0,0)}]
    \filldraw[black] (0,0) circle (\ptr pt);
    \filldraw[black] (1,0) circle (\ptr pt);
    \filldraw[black] (1.62,0) circle (\ptr pt);
    \filldraw[black] (2.62,0) circle (\ptr pt);
    
    \node[anchor = north] at (0,0-0.5*\h) {$p$};
    \node[anchor = north] at (1,0-0.5*\h) {$q$};
    \node[anchor = north] at (1.62,0-0.5*\h) {$r$};
    \node[anchor = north] at (2.62,0-0.5*\h) {$s$};
    
    \glebarc[0,1.62]
    \glebarc[1,1.62]
    
    \node[anchor = north] at (1.31,0-2*\h) {$(p<q<r<s)$};
    \end{scope}
    
    \begin{scope}[shift={(+5,0)}]
    \filldraw[black] (0,0) circle (\ptr pt);
    \filldraw[black] (1,0) circle (\ptr pt);
    \filldraw[black] (2,0) circle (\ptr pt);
    \filldraw[black] (2.62,0) circle (\ptr pt);
    
    \node[anchor = north] at (0,0-0.5*\h) {$p$};
    \node[anchor = north] at (1,0-0.5*\h) {$q$};
    \node[anchor = north] at (2,0-0.5*\h) {$r$};
    \node[anchor = north] at (2.62,0-0.5*\h) {$s$};
    
    \glebarc[0,2.62]
    \glebarc[1,1]
    
    \node[anchor = north] at (1.31,0-2*\h) {$(p<q<r<s)$};
    \end{scope}
    \end{tikzpicture}
\end{center}
\caption{A pair of consecutive curves $C_j$ and $C_{j+1}$ in the $C$-sequence are in one of the three configurations. The configuration on the left corresponds to outcomes (i) and (ii) in Lemma~\ref{l:consecutive-curves}. The middle one corresponds to outcomes (iii) and (iv). The right one corresponds to outcomes (v) and (vi).}\label{fig:consec-allow}
\end{figure}

To visualize Lemma~\ref{l:consecutive-curves}, consider Figure~\ref{fig:consec-allow}. If $C_j$ and $C_{j+1}$ are curves from the $C$-sequence, then they must follow one of the three configurations shown in the figure. (For each configuration, there are two choices for which curve is $C_j$ and which curve is $C_{j+1}$, so there are six outcomes in total as listed in Lemma~\ref{l:consecutive-curves}.)

\begin{figure}
\begin{center}
\def\ptr{2}
\def\h{0.5}
\def\glebarc[#1,#2]{
\draw[line width = 1 pt] (#1,0) arc(180:0:{0.5*#2} and {0.5*#2});}
    \begin{tikzpicture}
    \begin{scope}[shift={(-3,+2)}]
    \filldraw[black] (0,0) circle (\ptr pt);
    \filldraw[black] (1.62,0) circle (\ptr pt);
    \filldraw[black] (2.62,0) circle (\ptr pt);
    
    \node[anchor = north] at (0,0-0.5*\h) {$p$};
    \node[anchor = north] at (1.62,0-0.5*\h) {$q$};
    \node[anchor = north] at (2.62,0-0.5*\h) {$r$};
    
    \glebarc[0,1.62]
    \glebarc[0,2.62]
    
    \node[anchor = north] at (1.31,0-2*\h) {$(p<q<r)$};
    \end{scope}
    
    \begin{scope}[shift={(+3,+2)}]
    \filldraw[black] (0,0) circle (\ptr pt);
    \filldraw[black] (1.62,0) circle (\ptr pt);
    \filldraw[black] (2.62,0) circle (\ptr pt);
    
    \node[anchor = north] at (0,0-0.5*\h) {$p$};
    \node[anchor = north] at (1.62,0-0.5*\h) {$q$};
    \node[anchor = north] at (2.62,0-0.5*\h) {$r$};
    
    \glebarc[0,2.62]
    \glebarc[1.62,1]
    
    \node[anchor = north] at (1.31,0-2*\h) {$(p<q<r)$};
    \end{scope}
    
    \begin{scope}[shift={(-3,-2)}]
    \filldraw[black] (0,0) circle (\ptr pt);
    \filldraw[black] (2.62,0) circle (\ptr pt);
    
    \node[anchor = north] at (0,0-0.5*\h) {$p$};
    \node[anchor = north] at (2.62,0-0.5*\h) {$q$};
    
    \glebarc[0,2.62]
    \draw[line width = 1 pt] (0,0) arc(180:0:{0.5*2.62} and {0.35*2.62});
    
    \node[anchor = north] at (1.31,0-2*\h) {$(p<q)$};
    \end{scope}
    
    \begin{scope}[shift={(+3,-2)}]
    \filldraw[black] (0,0) circle (\ptr pt);
    \filldraw[black] (1,0) circle (\ptr pt);
    \filldraw[black] (2,0) circle (\ptr pt);
    \filldraw[black] (2.62,0) circle (\ptr pt);
    
    \node[anchor = north] at (0,0-0.5*\h) {$p$};
    \node[anchor = north] at (1,0-0.5*\h) {$q$};
    \node[anchor = north] at (2,0-0.5*\h) {$r$};
    \node[anchor = north] at (2.62,0-0.5*\h) {$s$};
    
    \glebarc[0,1]
    \glebarc[2,0.62]
    
    \node[anchor = north] at (1.31,0-2*\h) {$(p<q<r<s)$};
    \end{scope}

    \end{tikzpicture}
\end{center}
\caption{If $C_j$ and $C_{j+1}$ are curves from the $C$-sequence, then they cannot follow any of these four forbidden configurations.}\label{fig:consec:forbid}
\end{figure}
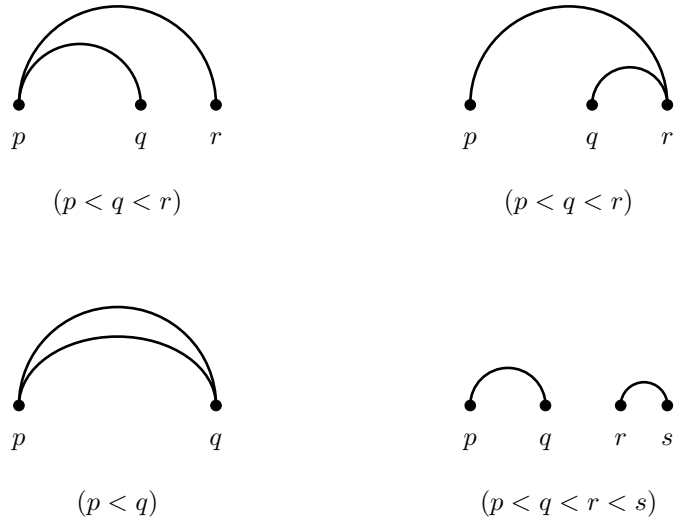

Equivalently, Lemma~\ref{l:consecutive-curves} states that the four configurations shown in Figure~\ref{fig:consec:forbid} cannot represent two consecutive curves in the $C$-sequence.

\begin{proof}[Proof of Lemma~\ref{l:consecutive-curves}]
Suppose, for the sake of contradiction, that none of the six outcomes happens. Then either $a=c$, or $b=d$, or $b < c$, or $d < a$. Recall that we obtain $A_j$ from $A_{j-1}$ by removing $a$ and adding $b$, and we obtain $A_{j+1}$ from $A_j$ by removing $c$ and adding $d$. If the first case, $a = c$, happens, then we would need two copies of $a$ in the set $A_{j-1}$, a contradiction. Similarly, if the second case, $b = d$, happens, then we would need two copies of $b$ in the set $A_{j+1}$, a contradiction. If the third case, $b < c$, or the fourth case, $d < a$, happens, then there would be an element $x \in [n]$, different from $a, b, c, d$, such that $A_{j-1} = \{x, a, c\}$ and $A_{j+1} = \{x, b, d\}$. Note that $A_{j-1}$ and $A_{j+1}$ are not weakly separated, a contradiction.
\end{proof}

\begin{proof}[Proof of Proposition~\ref{prop:simple-BiGP}]
For any pair of different nodes $i, j \in [n]$, we know that the pair is connected by at most one black edge. It suffices to show that for each new red edge added with endpoints $i$ and $j$, the nodes $i$ and $j$ have not already had a black edge or a red edge connecting them.

Recall that the red edges are added in $N-1$ different steps. Consider the red edges added in the $t^{\text{th}}$ step. Following Definition~\ref{defn:BiGP}, we consider the pairs $C_t = \{a,b\}$ and $C_{t+1} = \{c,d\}$, with $a < b$ and $c < d$. Let $\chi_t$ and $\chi_{t+1}$ denote the crossings in the wiring diagram which correspond to $C_t$ and $C_{t+1}$, respectively. We claim that the red edges constructed in this step correspond to crossings which happen between $\chi_t$ and $\chi_{t+1}$ (on different levels: $s_{k'}$ with $k' \neq 3$). With this claim, the proposition is proved, because we are selecting different crossings in each of the $N-1$ steps.

To establish the claim, we use Lemma~\ref{l:consecutive-curves}. The pairs $C_t$ and $C_{t+1}$ exhibit one of the six outcomes as listed in the lemma. Consider the outcome (iii) (and one argues similarly for the other outcomes). Note that wires $a$ and $b$ cross at $\chi_t$, wires $c$ and $d$ cross at $\chi_{t+1}$, and no other crossings can happen on the third level. This means that wires $a$ and $d$ must cross somewhere between $\chi_t$ and $\chi_{t+1}$ on the $(k')^{\text{th}}$ level for some $k' \ge 4$. Similarly, wires $b$ and $c$ must cross somewhere between $\chi_t$ and $\chi_{t+1}$ as well on either the first or the second level. 
\end{proof}

The decomposition algorithm (Algorithm~\ref{alg:fish}) uses the result that in $\mathbf{D}(P)$ the length of each $\mu(H_i)$ is at most $4$, which follows from the following proposition.

\begin{figure}
\begin{center}
\def\ptr{2}
\def\h{0.35}
\def\glebarc[#1,#2]{
\draw[line width = 1 pt] (#1,0) arc(180:0:{0.5*#2} and {0.5*#2});}
    \begin{tikzpicture}[scale = 1.6]
    
    \filldraw[black] (0,0) circle (\ptr pt);
    \filldraw[black] (1,0) circle (\ptr pt);
    \filldraw[black] (2,0) circle (\ptr pt);
    \filldraw[black] (3,0) circle (\ptr pt);
    \filldraw[black] (4,0) circle (\ptr pt);
    \filldraw[black] (5,0) circle (\ptr pt);
    
    \node at (0,0-\h) {$i$};
    \node at (1,0-\h) {$i+1$};
    \node at (2,0-\h) {$i+2$};
    \node at (3,0-\h) {$i+3$};
    \node at (4,0-\h) {$i+4$};
    \node at (5,0-\h) {$i+5$};
    
    \glebarc[0,1]
    \glebarc[1,1]
    \glebarc[2,1]
    \glebarc[3,1]
    \glebarc[4,1]
    \glebarc[0,2]
    \glebarc[1,2]
    \glebarc[2,2]
    \glebarc[3,2]
    
    \end{tikzpicture}
\end{center}
\caption{The nine curves in this configuration cannot simultaneously appear in the \Gleb diagram $\mathbf{D}(P)$ of a monotone weakly separated path $P$.}\label{fig:nine-curve-lemma}
\end{figure}
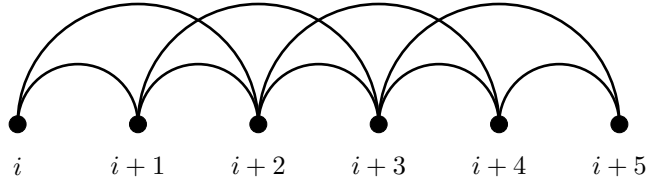

\begin{proposition}\label{prop:mu(H_i)-le-4}
In $\mathbf{D}(P)$, there is no index $i$ such that all nine (black) edges $\{i,i+1\}$, $\{i+1, i+2\}$, $\{i+2,i+3\}$, $\{i+3,i+4\}$, $\{i+4,i+5\}$, $\{i,i+2\}$, $\{i+1,i+3\}$, $\{i+2,i+4\}$, and $\{i+3,i+5\}$ appear. (See Figure~\ref{fig:nine-curve-lemma} for an illustration.)
\end{proposition}

Before proceeding to the proof, we give a quick explanation here how this proposition implies that each $\mu(H_i)$ in Algorithm~\ref{alg:fish} is at most $4$. If $\mu(H_i)$ is at least $5$, then there must be an index $i$ for which the five intervals $[i,i+1], [i+1,i+2], \ldots, [i+4,i+5]$ exceed their weight limits. It is not hard to see that this implies $\wt([i,i+1]) = \cdots = \wt([i+4,i+5]) = 2$, and thus arcs $\{i+j-1,i+j\}$ exist in $\mathbf{D}(P)$ for $j \in [5]$, and arcs $\{i+j-2,i+j\}$ exist in $\mathbf{D}(P)$ for $j \in [6]$. These many arcs would contain the configuration as shown in Figure~\ref{fig:nine-curve-lemma}.

\begin{proof}[Proof of Proposition~\ref{prop:mu(H_i)-le-4}]
Suppose, for the sake of contradiction, that there is such an index $i$. Since there are at most three pieces of curves above each unit interval, we know that there are no more black edges above the segment $[i+1,i+4]$. Above $[i,i+1]$, we now have two black curves. Thus, there can be at most one more black curve whose right endpoint is $i+1$. Call this curve, if it exists, $\zeta$. Similarly, there is at most one curve connecting $i+4$ and some $j > i+5$. Call this curve, if it exists, $\mu$.

Call the nine curves in the proposition $\alpha$, $\beta$, $\gamma$, $\delta$, $\varepsilon$, $\eta$, $\theta$, $\kappa$, $\lambda$, in the same order as displayed in the proposition statement. These nine curves, together with $\zeta$ and $\mu$, are all black curves above the segment $[i,i+5]$. Recall that we have the $C$-sequence $C_1, C_2, \ldots, C_N$ which lists all the black curves in $\mathbf{BiD}(P)$. We will consider which two curves are consecutive in this sequence.

Consider the curve $\beta$. From Lemma~\ref{l:consecutive-curves}, we know that none of $\delta$, $\varepsilon$, $\eta$, $\theta$, $\lambda$, and ($\mu$) can be consecutive to $\beta$ in the $C$-sequence. (The parentheses about $\mu$ in the previous sentence serve as a reminder that perhaps $\mu$ does not exist.) Moreover, $\alpha$ and $\beta$ cannot be consecutive edges in the $C$-sequence. Otherwise, there would be a red edge connecting $i$ and $i+2$ in $\mathbf{BiD}(P)$, contradicting Proposition~\ref{prop:simple-BiGP} as the black curve $\eta$ is already connecting $i$ and $i+2$. Similarly, $\gamma$ cannot be a neighbor of $\beta$.

There are only two choices left for the neighbors of $\beta$: $(\zeta)$ and $\kappa$. If $\zeta$ does exist, then $\beta$ cannot be $C_1$ (the starting curve in the $C$-sequence). We know $\beta$ cannot be $C_N$ either. Thus, $\beta$ must be adjacent to {\em both} $\zeta$ and $\kappa$. If $\zeta$ does not exist, $\beta$ must be adjacent to $\kappa$. In either case, we know $\beta$ and $\kappa$ are neighbors in the $C$-sequence, and thus there must be a red curve connecting $i+1$ and $i+4$ corresponding to a crossing between $\beta$ and $\kappa$.

However, the same reasoning implies that $\theta$ and $\delta$ must be adjacent in the $C$-sequence as well. There must be another red curve connecting $i+1$ and $i+4$ corresponding to a crossing between $\theta$ and $\delta$. This contradicts Proposition~\ref{prop:simple-BiGP}.
\end{proof}

\begin{lem}\label{l:mu(H_i)=1}
Suppose that $i$ is an integer with $4 \le i \le n-4$. Suppose that in $\mathbf{D}(P)$, there are indices $i'$ and $i''$ with $i' \le i - 3$ and $i'' \ge i + 4$ such that the edges $\{i-1,i\}$, $\{i,i+1\}$, $\{i+1,i+2\}$, $\{i-1,i+1\}$, $\{i,i+2\}$, $\{i',i\}$, and $\{i+1, i''\}$ appear in the diagram. Let $\alpha$ denote the curve connecting $i-1$ and $i$. Let $\delta$ denote the curve connecting $i+1$ and $i+2$. Let $(\beta)$ and $\gamma$ be the curves whose right endpoints are $i-1$. Let $(\varepsilon)$ and $\eta$ be the curves whose left endpoints are $i+2$.

Then either
\begin{itemize}
\item the neighbors of $\alpha$ in the $C$-sequence are $(\beta)$ and $\gamma$, or
\item the neighbors of $\delta$ in the $C$-sequence are $(\varepsilon)$ and $\eta$.
\end{itemize}
\end{lem}

Once again, the parentheses about $\beta$ and $\varepsilon$ in the lemma above mean ``if it exists". In the case $i = 4$, there is only one curve whose right endpoint is $i-1 = 3$. We denote that curve by $\gamma$, and $\beta$ is non-existent. Similarly, $\varepsilon$ is non-existent if and only if $i = n-4$.

\begin{proof}[Proof of Lemma~\ref{l:mu(H_i)=1}]
Let $\xi$ denote the curve connecting $i-1$ and $i+1$, and let $\zeta$ denote the curve connecting $i$ and $i+2$. Suppose that the neighbors of $\alpha$ are not $(\beta)$ and $\gamma$. Then by \Gleb diagram chasing, $\zeta$ must be a neighbor of $\alpha$. Therefore, in $\mathbf{BiD}(P)$, we have a red curve connecting $i-1$ and $i+2$ corresponding to a crossing between the crossings of $\alpha$ and $\zeta$. This shows that $\xi$ cannot be a neighbor of $\delta$. Thus, the neighbors of $\delta$ are $(\varepsilon)$ and $\eta$.
\end{proof}

The following lemma is a degenerate version of Lemma~\ref{l:mu(H_i)=1}. The proof is essentially the same as that of the previous lemma, so we omit it.

\begin{lem}\label{l:mu(H_i)=1-deg}
We have the following properties of $\mathbf{D}(P)$.
\begin{itemize}
\item[(a)] Suppose that $n \ge 7$. Suppose that there is an index $i \ge 7$ such that the edges $\{1,2\}$, $\{2,3\}$, $\{3,4\}$, $\{4,5\}$, $\{2,4\}$, $\{3,5\}$, and $\{4,i\}$ appear in $\mathbf{D}(P)$. Let $\alpha$ denote the curve connecting $4$ and $5$. Let the black curves whose left endpoints are $5$ be $(\beta)$ and $\gamma$. Then the neighbors of $\alpha$ are $(\beta)$ and $\gamma$.
\item[(b)] Suppose that $n \ge 7$. Suppose that there is an index $i \le n-6$ such that the edges $\{n,n-1\}$, $\{n-1,n-2\}$, $\{n-2,n-3\}$, $\{n-3,n-4\}$, $\{n-1,n-3\}$, $\{n-2,n-4\}$, and $\{n-3,i\}$ appear in $\mathbf{D}(P)$. Let $\alpha$ denote the curve connecting $n-3$ and $n-4$. Let the black curves whose right endpoints are $n-4$ be $(\beta)$ and $\gamma$. Then the neighbors of $\alpha$ are $(\beta)$ and $\gamma$.
\item[(c)] When $n=6$, the edges $\{1,2\}$, $\{2,3\}$, $\{3,4\}$, $\{4,5\}$, $\{5,6\}$, $\{2,4\}$, $\{3,5\}$ cannot simultaneously appear in $\mathbf{D}(P)$.
\end{itemize}
\end{lem}

Lemmas~\ref{l:mu(H_i)=1} and \ref{l:mu(H_i)=1-deg} deal with the situation where we encounter one unit interval of weight $2$. When there are two consecutive unit intervals of weight $2$, \Gleb diagram chasing gives a result similar to Lemma~\ref{l:mu(H_i)=1} as follows.

\begin{lem}\label{l:mu(H_i)=2}
Suppose that $i$ is an integer such that $4 \le i \le n-5$. Suppose that in $\mathbf{D}(P)$, there are indices $i'$ and $i''$ with $i' \le i - 3$ and $i'' \ge i+5$ such that the edges $\{i-1,i\}$, $\{i,i+1\}$, $\{i+1,i+2\}$, $\{i+2,i+3\}$, $\{i-1,i+1\}$, $\{i,i+2\}$, $\{i+1,i+3\}$, $\{i',i\}$, and $\{i+2,i''\}$ appear in the diagram. Let $\alpha$ denote the curve connecting $i+2$ and $i+3$. Let the curves whose left endpoints are $i+3$ be $(\beta)$ and $\gamma$. Let $\delta$ denote the curve connecting $i-1$ and $i$. Let the curves whose right endpoints are $i-1$ be $(\varepsilon)$ and $\eta$. 

Then both of the following are true:
\begin{itemize}
\item the neighbors of $\alpha$ are $(\beta)$ and $\gamma$.
\item the neighbors of $\delta$ are $(\varepsilon)$ and $\eta$.
\end{itemize}
\end{lem}

A degenerate version of Lemma~\ref{l:mu(H_i)=2} is Lemma~\ref{l:mu(H_i)=2-deg} below.

\begin{lem}\label{l:mu(H_i)=2-deg}
We have the following properties of $\mathbf{D}(P)$.
\begin{itemize}
\item[(a)] The edges $\{1,2\}$, $\{2,3\}$, $\{3,4\}$, $\{4,5\}$, $\{5,6\}$, $\{2,4\}$, $\{3,5\}$, and $\{4,6\}$ cannot simultaneously appear in $\mathbf{D}(P)$.
\item[(b)] The edges $\{n-5,n-4\}$, $\{n-4,n-3\}$, $\{n-3,n-2\}$, $\{n-2,n-1\}$, $\{n-1,n\}$, $\{n-5,n-3\}$, $\{n-4,n-2\}$, and $\{n-3,n-1\}$ cannot simultaneously appear in $\mathbf{D}(P)$.
\end{itemize}
\end{lem}

Given an \Gleb diagram $\mathbf{D}(P)$, we have seen in Algorithm~\ref{alg:fish} that we can decompose the interval $[1,n]$ into
\[
[1,n] = L_0 \cup H_1 \cup L_1 \cup \cdots \cup L_m,
\]
where every unit interval in $L_i$ does not exceed its weight limit, and every unit interval in $H_i$ has weight $2$. (In Algorithm~\ref{alg:fish}, we defined this decomposition for \Gleb diagrams with at least one unit interval with weight $2$. Here, we define it for any $\mathbf{D}(P)$. For \Gleb diagrams in which every unit interval does not exceed its weight limit, we can simply let $m = 0$ and $L_0 = L_m = [1,n]$.) The following proposition gives some restrictions on the lengths of the intervals $L_0, H_1, L_1, \ldots, H_m, L_m$. 

\begin{proposition}\label{prop:lengths}
We have
\begin{itemize}
\item[(a)] $\mu(L_0) \ge 2$,
\item[(b)] $\mu(L_m) \ge 2$,
\item[(c)] for $1 \le i \le m-1$, the interval $L_i$ satisfies $\mu(L_i) \ge 3$, and
\item[(d)] for $1 \le i \le m$, the interval $H_i$ satisfies $\mu(H_i) \le 4$.
\end{itemize}
\end{proposition}

\begin{proof}
\textbf{(a)} and \textbf{(b)} are clear, because the unit intervals $[1,2]$, $[2,3]$, $[n-2,n-1]$, $[n-1,n]$ never exceed their weight limits, by definition of $\wtlim$.

\textbf{(c)}. We will show that $\mu(L_i)$ cannot be $1$ or $2$. First, suppose $\mu(L_i) = 1$. Then there is some index $j$ such that $L_i = [j+2,j+3]$. Since $\wt([j+1,j+2]) = \wt([j+3,j+4]) = 2$, the edges $\{j+1,j+2\}$, $\{j+3,j+4\}$, $\{j,j+2\}$, $\{j+1,j+3\}$, $\{j+2,j+4\}$, and $\{j+3,j+5\}$ must appear in $\mathbf{D}(P)$. After drawing these six curves, we see that there are now three pieces of curves above $[j+1,j+2]$ and also there are now three pieces of curves above $[j+3,j+4]$. At the moment, there are only two pieces of curves above $[j+2,j+3]$, and thus there must be another piece of curve above $[j+2,j+3]$. Since there can be no more curves above $[j+1,j+2] \cup [j+3,j+4]$, the only option is to connect $j+2$ and $j+3$. However, this would make $\wt([j+2,j+3]) = 2$, a contradiction.

Second, suppose $\mu(L_i) = 2$. Then there is some index $j$ such that $L_i = [j+2,j+4]$. Note that $\wt([j+1,j+2]) = \wt([j+4,j+5]) = 2$. By a similar argument as in the previous case, we know that the following three pairs $\{j+2, j+3\}$, $\{j+3, j+4\}$, and $\{j+2,j+4\}$ must be connected by edges. However, this would make $\wt([j+2,j+3]) = \wt([j+3,j+4]) = 2$, a contradiction.

\textbf{(d)} follows from Proposition~\ref{prop:mu(H_i)-le-4}.
\end{proof}

\begin{proof}[Proof of Theorem~\ref{thm:alg-good}(a)]
In Algorithm~\ref{alg:fish}, we note that each interval we add to $\mathcal{F}(\mathbf{D})$ contains either one or two unit intervals from $\bigcup_{i=0}^m L_i$. More precisely, Case~3.3 and Case~4.3 in the algorithm are the only two cases that give intervals with two unit intervals from $\bigcup_{i=0}^m L_i$. Let $F$ and $F'$ be two different intervals in $\mathcal{F}(\mathbf{D})$. From Proposition~\ref{prop:lengths}(c), we see that if either $F$ or $F'$ does not come from these two cases, then $\mu(F \cap F') = 0$. The only potentially problematic case is when both $F$ and $F'$ come from Case~3.3 or Case~4.3 and the overlap $F \cap F'$ has length $1$. We will show that this is not possible.

Suppose, for the sake of contradiction, that $F$ and $F'$ share an interior point. Then there must be an index $a$ such that $F \cap F' = [a+3, a+4]$. This means $\wt([a+1,a+2]) = \wt([a+5,a+6]) = 2$, and each of the three unit intervals $[a+2,a+3]$, $[a+3,a+4]$, and $[a+4,a+5]$ has weight at most $11/6$. Since both $F$ and $F'$ come from Case~3.3 or Case~4.3, we see that the following edges $\{a+1,a+2\}$, $\{a+2,a+3\}$, $\{a+4,a+5\}$, $\{a+5,a+6\}$, $\{a,a+2\}$, $\{a+1,a+3\}$, $\{a+4,a+6\}$, and $\{a+5,a+7\}$ appear in $\mathbf{D}(P)$.

After drawing these eight edges, we observe that there are already three pieces of curve above $[a+1,a+2]$ and another three pieces above $[a+5,a+6]$. At the moment, there are only two pieces of curve above $[a+2,a+3]$. As no more curve can be added above $[a+1,a+2]$, there must be another curve $\gamma$ whose left endpoint is $a+2$. The right endpoint must be either $a+4$ or $a+5$. However, if the right endpoint were $a+4$, then $\wt([a+2,a+3])$ would be $2$, a contradiction. This forces $\gamma$ to connect $a+2$ and $a+5$. Now, there are three pieces of curves above $[a+2,a+3]$, and also three pieces above $[a+4,a+5]$.

Now consider $[a+3,a+4]$. At the moment, there is only one piece of curve above it, and so we need two more pieces. On the other hand, no more curves can be added above $[a+1,a+3]$ or above $[a+4,a+6]$. This gives a contradiction. We have finished the proof of Theorem~\ref{thm:alg-good}(a).
\end{proof}

We have shown that the intervals in the collection $\mathcal{F}(\mathbf{D})$ do not overlap. Next, we show that each interval has weight under or equal to its weight limit. 

\begin{proof}[Proof of Theorem~\ref{thm:alg-good}(b)]
Let $I$ be an arbitrary interval in $\mathcal{F}(\mathbf{D})$.

\smallskip

\underline{Case~1.} Suppose that $I$ comes from some $H_i = [a,a+4]$. We have that the edges $\{a,a+1\}$, $\{a+1, a+2\}$, $\{a+2, a+3\}$, $\{a+3, a+4\}$, $\{a-1,a+1\}$, $\{a,a+2\}$, $\{a+1, a+3\}$, $\{a+2, a+4\}$, and $\{a+3,a+5\}$ appear in $\mathbf{D}$. Since there are three pieces of curve above $[a,a+1]$, we know that $a \ge 3$. With Proposition~\ref{prop:mu(H_i)-le-4}, we know that there is no edge connecting $a-1$ and $a$ in $\mathbf{D}$. This means that if $a \ge 4$, we have
\[
\wt([a-1,a]) \le \frac{1}{2} + \frac{1}{2} + \frac{1}{3} = \frac{4}{3},
\]
whence $\wt([a-1,a+2]) \le \frac{4}{3} + 2 + 2 = \frac{16}{3} < \frac{11}{2} = \wtlim([a-1,a+2])$. If $a = 3$, we have $\wt([a-1,a]) \le \frac{1}{2} + \frac{1}{2} = 1$, whence $\wt([a-1,a+2]) \le 1 + 2 + 2 = 5 < \frac{31}{6} = \wtlim([a-1,a+2])$. Thus, if $I = [a-1,a+2]$, we have shown that $\wt(I) \le \wtlim(I)$. On the other hand, if $I = [a+2,a+5]$, the argument is analogous.

\smallskip

\underline{Case~2.1.} We have some index $a$ such that $H_i = [a,a+3]$ and $I = [a-1,a+3]$. There is no edge connecting $a-1$ and $a$. Since there are three pieces of curve above $[a,a+1]$, we have $a \ge 3$. If $a \ge 4$, then
\[
\wt([a-1,a+3]) \le \left( \frac{1}{2} + \frac{1}{2} + \frac{1}{3} \right) + 2 + 2 + 2 = \frac{22}{3} = \wtlim([a-1,a+3]).
\]
If $a = 3$, then $\wt([a-1,a+3]) \le \left( \frac{1}{2} + \frac{1}{2} \right) + 2 + 2 + 2 = 7 = \wtlim([a-1,a+3])$.

\smallskip

\underline{Case~2.2.} In this case, $H_i = [a,a+3]$ and $I = [a,a+4]$. There is an edge connecting $a-1$ and $a$. Therefore, by Proposition~\ref{prop:mu(H_i)-le-4}, there is no edge connecting $a+3$ and $a+4$. Hence, the weight calculation is similar to Case~2.1.

\smallskip

\underline{Case~3.1} and \underline{Case~3.2} are also similar.

\smallskip

\underline{Case~3.3.} In this case, $H_i = [a,a+2]$. Let us show that $\wt(I) \le \wtlim(I)$ for $I = [a-2,a+1]$. By symmetry, the case when $I = [a+1,a+4]$ is analogous. The following edges $\{a-1,a\}$, $\{a,a+1\}$, $\{a+1,a+2\}$, $\{a+2,a+3\}$, $\{a-1,a+1\}$, $\{a,a+2\}$, and $\{a+1,a+3\}$ appear in $\mathbf{D}(P)$. Since there are three pieces of curve above $[a,a+1]$, we have $a \ge 3$. If $a = 3$, we have no more curves above $[a-1,a] = [2,3]$, and so there must be a curve connecting $a-2 = 1$ and $a-1 = 2$, contradicting Lemma~\ref{l:mu(H_i)=2-deg}(a). Therefore, $a \ge 4$. There must be one more curve whose right endpoint is $a$. The left endpoint cannot be $a-2$; otherwise $\wt([a-1,a])$ would be $2$. Thus, there is an index $a' \le a-3$ such that there is a curve connecting $a'$ and $a$.

Similarly, we find that $a \le n-5$ and there is an index $a'' \ge a+5$ such that there is a curve connecting $a+2$ and $a''$. Note that we now have the assumptions of Lemma~\ref{l:mu(H_i)=2} (with $a, a', a''$ here playing the roles of $i, i', i''$ in the lemma). Following the notations in the lemma, let $\delta$ denote the curve connecting $a-1$ and $a$. Let the curves whose right endpoints are $a-1$ be $(\varepsilon)$ and $\eta$. By Lemma~\ref{l:mu(H_i)=2}, we have that the neighbors of $\delta$ are $(\varepsilon)$ and $\eta$.

Let $\kappa$ denote the curve connecting $a'$ and $a$. We have that the left endpoints of $\kappa$, $(\varepsilon)$, and $\eta$ are all distinct. (Otherwise, by a little bit of \Gleb diagram chasing, there would be a red curve with the same endpoints as $\kappa$ in $\mathbf{BiD}(P)$, contradicting Proposition~\ref{prop:simple-BiGP}.) We now check the weight $\wt(I)$.

If $a = 4$, the curve $\varepsilon$ is non-existent. The curve $\kappa$ connects $1$ and $4$. The curve $\eta$ connects $2$ and $3$. We have $\wt(I) = \wt([2,5]) = \frac{31}{6} = \wtlim([2,5])$.

If $a \ge 5$, then $\varepsilon$ exists. Suppose that the lengths of $\kappa$, $\varepsilon$, $\eta$ are $u+1$, $v$, $w$, respectively. Since the left endpoints of the three curves are all distinct, we have that $u$, $v$, $w$ are distinct positive integers. It is straightforward to compute $\wt(I) = \frac{2}{u+1} + \frac{1}{v} + \frac{1}{w} + \frac{7}{2}$. Note that the weight limit is $\wtlim(I) = \frac{11}{2}$. It is a pleasant exercise to show that for distinct positive integers $u, v, w$, we have the inequality
\[
\frac{2}{u+1} + \frac{1}{v} + \frac{1}{w} \le 2,
\]
which we will leave to the reader. This shows that $\wt(I) \le \wtlim(I)$.

\smallskip

\underline{Case~4.1} and \underline{Case~4.2} are also similar to Case~2.1 above.

\smallskip

\underline{Case~4.3.} In this case, $H_i = [a,a+1]$. We would like to show that either $\wt([a-2,a+1]) \le \wtlim([a-2,a+1])$ or $\wt([a,a+3]) \le \wtlim([a,a+3])$. Since $\wt([a,a+1]) = 2$, we have that $3 \le a \le n-3$. Start by considering edge cases. If $n = 6$, then $a = 3$ and the edges $\{1,2\}$, $\{2,3\}$, $\{3,4\}$, $\{4,5\}$, $\{5,6\}$, $\{2,4\}$, and $\{3,5\}$ appear in $\mathbf{D}(P)$. This directly contradicts Lemma~\ref{l:mu(H_i)=1-deg}(c). Assume now that $n \ge 7$. If $a=3$, then we are in the situation of Lemma~\ref{l:mu(H_i)=1-deg}(a). By using an argument similar to one in Case~3.3, we find that $[3,6]$ is under its weight limit. The case $a = n-3$ is analogous.

Now assume $4 \le a \le n-4$. We see that there exist indices $a' \le a-3$ and $a'' \ge a+4$ such that there are edges $\{a',a\}$ and $\{a+1,a''\}$ in $\mathbf{D}(P)$. We are in the situation of Lemma~\ref{l:mu(H_i)=1} (with $a, a', a''$ here playing the roles of $i, i', i''$ in the lemma). Following notations in the lemma, let us denote the curve connecting $a-1$ and $a$ by $\alpha$, and denote the curve connecting $a+1$ and $a+2$ by $\delta$. Let the curves whose right endpoints are $a-1$ be $(\beta)$ and $\gamma$. Let the curves whose left endpoints are $a+2$ be $(\varepsilon)$ and $\eta$. Lemma~\ref{l:mu(H_i)=1} says that either the neighbors of $\alpha$ are $(\beta)$ and $\gamma$, or the neighbors of $\delta$ are $(\varepsilon)$ and $\eta$.

If the neighbors of $\alpha$ are $(\beta)$ and $\gamma$, then by an argument similar to one in Case~3.3, we find that $[a-2,a+1]$ has weight under or equal to its weight limit. On the other hand, if the neighbors of $\delta$ are $(\varepsilon)$ and $\eta$, then the weight of $[a,a+3]$ is under or equal to its weight limit. We have finished the proof.
\end{proof}

\subsection{Repeatable patterns for $k=3$}

We now establish the lower bound on $\Ans(3,n)$ by giving explicit repeatable patterns.

\begin{defn}
Let $P = (A_0, A_1, \ldots, A_N)$ and $Q = (B_0, B_1, \ldots, B_M)$ be sequences of $k$-element sets of integers. Suppose that there exists an integer $t$ such that $A_N = B_0 + t$. Then we define the {\em concatenation} of $P$ and $Q$ to be the sequence
\[
P \ast Q := \left( A_0, A_1, \ldots, A_{N-1}, B_0 + t, B_1 + t, \ldots, B_M + t \right).
\]
\end{defn}

Therefore, a repeatable pattern $R$ is a monotone weakly separated path such that for any positive integer $m$, the $m^{\text{th}}$-concatenation power $R \ast R \ast \cdots \ast R$ of $R$ is well-defined and is also a monotone weakly separated path. 

Now we construct optimal monotone weakly separated paths as follows. We define:
\begin{align*}
&P_4 = 123-124-134-234, \\
&P_5 = 123-124-125-145-245-345, \\
&P_6 = 123-124-125-145-245-345-456, \\
&P_7 = 123-124-125-145-245-345-456-457-567, \\
&P_8 = 123-124-125-145-245-345-456-457-567-578-678, \text{ and} \\
&P_9 = \\
&123-124-125-145-245-345-456-457-567-578-579-589-789.
\end{align*}
We also define
\[
P = 123-124-125-145-245-345-456-457-567-578-678-789.
\]
It is straightforward to check that $P$ is a repeatable pattern with parameters $(L,d) = (11,6)$.

For each integer $n \ge 10$, define $P_n := P \ast P_{n-6}$. It is also straightforward to check that for every integer $n \in \mathbb{Z}_{\ge 4}$, the sequence $P_n$ is a monotone weakly separated path from $\{1, 2, 3\}$ to $\{n-2, n-1, n\}$ with $\left\lceil \frac{11}{6} n \right\rceil - 4$ terms. Combining these constructions with Corollary~\ref{cor:M3n-upper}, we have proved the following theorem.

\begin{thm}\label{thm:M3n}
For each positive integer $n \ge 4$, we have
\[
\mathcal{M}(3,n) = \left\lceil \frac{11}{6} n \right\rceil - 5.
\]
\end{thm}

We end this section with a remark about the general formula for $\mathcal{M}(k,n)$. Considering the formulas for $k = 1, 2, 3$, one might conjecture that in general there exist real numbers $a_k$ and $b_k$ for which the formula $\mathcal{M}(k,n) = \left\lceil a_k n + b_k \right\rceil$ holds for every $n \ge k+1$. Unfortunately, from our computational results, we can show, for example by Fourier--Motzkin Elimination, that there cannot be such a formula when $k = 4$. The formula for $\mathcal{M}(4,n)$ has to be somewhat more complicated.

\section{Periodicity for $\Ans(k,n)$}\label{sec:rationality}
Fix a positive integer $k$ throughout this section. Define the constant \[c_k := \lim_{n \to \infty} \frac{\Ans(k,n)}{n}.\] 
Our main goal of this section is to show that for any $k$, $c_k$ exists, is rational and can be achieved by repeatable patterns.

\subsection{Existence of $c_k$}
In Section~\ref{sec:bounds}, $c_k$ was explicitly computed for $k=1,2,3$. To be precise, $c_1 = 1$, $c_2 = \frac{3}{2}$, and $c_3 = \frac{11}{6}$. Theorem~\ref{thm:existence} shows that this limit exists for all $k \in \mathbb{N}$. 

\begin{thm}
\label{thm:existence}
The limit $c_k$ exists  
for any $k\in \mathbb{N}$.
\end{thm}

The proof rests on the following lemma.

\begin{lem} \label{lem:superadditivity}
For positive integers $k<n\leq m$, we have
\begin{itemize}
\item[(a)] $\Ans(k,n)\leq \Ans(k,m),$ and
\item[(b)] $\Ans(k,n)+ \Ans(k,m)\leq \Ans(k,n+m).$
\end{itemize}
\end{lem}
\begin{proof}
We prove part (b); part (a) follows since $\Ans(k,m)$ is nonnegative. It follows from the basic theory of Coxeter groups (see \cite{bjorner-brenti}) that the longest permutation $w_{0,n}$ in $S_n$ is less than the longest permutation $w_{0,n+m}$ in $S_{n+m}$ in the weak Bruhat order, if we view $S_n \subset S_{n+m}$ as those permutations fixing $n+1,\ldots,n+m$ pointwise. Thus we can write 
\[
w_{0,n+m} = w_{0,n} \cdot u
\]
with $\ell(w_{0,n+m})=\ell(w_{0,n})+\ell(u)$. The permutation $u=u_1 \ldots u_{n+m}$ has $u_1>\cdots>u_m$, so we may write $u=u' \cdot w_{0,m}$, again with $S_m \subset S_{n+m}$ embedded in the standard way, and again with lengths adding. For any reduced words $\mathbf{i}, \mathbf{i'}$ of $w_{0,n}$ and $w_{0,m}$ and a reduced word $\mathbf{j}$ for $u'$, this implies that the concatenation $\mathbf{i}\mathbf{j}\mathbf{i'}$ is a reduced word for $w_{0,n+m}$. In particular, by choosing $\mathbf{i},\mathbf{i'}$ to each maximize the occurrences of $s_k$, we have:
\[
\Ans(k,n)+ \Ans(k,m)\leq \Ans(k,n+m),
\]
for all $k$.
\end{proof}
\begin{proof}[Proof of Theorem~\ref{thm:existence}]
By Lemma~\ref{lem:superadditivity}(b), $\mathcal{M}(k,n)$ is subadditive function of $n$, so by Fekete's Lemma, the desired limit
\[
\lim_{n \to \infty} \frac{\Ans(k,n)}{n}
\]
exists.

%, we know that for any $k\leq n\leq m,$ we have $$\Ans(k,m)\geq \left\lfloor \frac{m}{n} \right\rfloor \Ans(k,n).$$
%Fix any $n$. Then $$\frac{\Ans(k,m)}{m}\geq \left\lfloor \frac{m}{n}\right\rfloor \frac{n}{m} \frac{\Ans(k,n)}{n}.$$ Since $\left\lfloor \frac{m}{n}\right\rfloor \frac{n}{m}$ tends to $1$ when $m$ goes to infinity, any accumulation point is at least $\frac{\Ans(k,n)}{n}$.
%By Corollary~\ref{cor:sqrt(2k)}, we know that $\frac{\Ans(k,n)}{n}<\sqrt{2k}$, i.e., the sequence $\frac{\Ans(k,n)}{n}$ is bounded. Hence, it has a limit.	
\end{proof}

\subsection{Asymptotic equivalence with the pseudoline $k$-set problem}
As discussed in Section~\ref{sec:intro-k-sets}, in the context of the ``$k$-set problem" it is natural to consider a related problem, namely the maximization of the total number of appearances of $s_k$ and $s_{n-k}$ in a reduced word. Let $\AnsD(k,n)$ be the maximal total number of appearances of $s_k$ and $s_{n-k}$ in the reduced words from $S_n$. The following theorem shows that the same slopes $c_k$ arise in this version.

\begin{thm} For any $k\in \mathbb{N}$, the following limit exists and is given by
\label{thm:dk}
$$\lim_{n\to +\infty}\frac{\AnsD(k,n)}{n}=\lim_{n\to +\infty}\frac{\Ans(k,n)}{n}=c_k.$$
\end{thm}
\begin{proof}
Consider any reduced word and its wiring diagram. We say that a wire has type $(i,j,\pm)$ if its highest position is $i$ and its lowest position is $j$, and $+\ (-)$ means that the highest position is to the left (right) of the lowest position.
Note that no two wires share the same type (otherwise they should intersect at least twice, but our word is reduced).
Let $a$ be the number of wires which were at some moment at one of the $k$ highest levels, and let $b$ be the number of wires which were at some moment at one of the $k$ lowest levels. At most $2k^2$ wires are counted by both $a$ and $b$, so $a+b\leq n+2k^2$.
Note that the number of $s_k$ depends only on these $a$ wires
and the number of $s_{n-k}$ depends only on these $b$ wires.
Hence, the number of appearances of $s_k$ and $s_{n-k}$ in this reduced word is at most $\Ans(k,a)+\Ans(k,b) \leq c_k a+c_k b\leq c_k(n+2k^2) $.

Therefore $\Ans(k,n)\leq \AnsD(k,n) \leq c_k(n+2k^2)$. Thus the limit $\lim_{n\to +\infty}\frac{\AnsD(k,n)}{n}$ exists and is equal to $c_k$.
\end{proof}

\begin{remark}
We can similarly define numbers $\Ans(S,n)$ and $\AnsD(S,n)$ for any finite subset $S\subset \mathbb{N}$ and $n\in \mathbb{N}$. Their asymptotics are still the same and well-defined, i.e., \[
\lim_{n\to +\infty}\frac{\Ans(S,n)}{n}=\lim_{n\to +\infty}\frac{\AnsD(S,n)}{n}\in \mathbb{R}.
\]
\end{remark}

\subsection{Generalized wiring diagrams}
We now work towards showing the rationality of $c_k$. We introduce the new tool of \textit{generalized wiring diagrams}; these are certain wiring diagrams with infinitely many wires which are sometimes allowed to ``go to infinity". Intuitively, these diagrams allow us to reason about wiring diagrams for all $n \in \mathbb{N}$ simultaneously.

\begin{defn}
A \textit{generalized wiring diagram} consists of countably many wires, labeled by $1,2,\ldots$ starting at \text{levels} $1,2,\ldots$ respectively from top to bottom, traveling from left to right such that at each timestamp $t$, either
\begin{itemize}
\item two wires at adjacent levels cross; or
\item one wire goes to infinity $\infty$, intersecting all wires at lower levels.
\end{itemize}
A generalized wiring diagram is \textit{reduced} if no pair of wires cross more than once.
\end{defn}
To clarify, when two wires at level $h$ and $h+1$ cross as in the usual wiring diagrams, we say that they cross at level $h$. And when a wire $a$ at level $h$ goes to infinity, we say that this wire \textit{falls}, and it creates intersections at levels $h,h+1,\ldots$, while the wires which were at levels $h+1,h+2,\ldots$ before wire $a$ falls go to levels $h,h+1,\ldots$ respectively, so that at every timestamp, there is a wire at each level indexed by positive integers. 

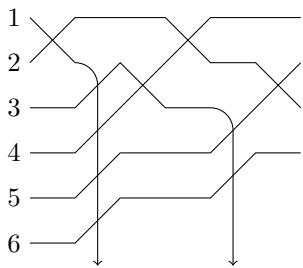
\begin{figure}[h!]
\centering
\begin{tikzpicture}[scale=0.6]
\node[left] at (0,0) {$1$};
\node[left] at (0,-1) {$2$};
\node[left] at (0,-2) {$3$};
\node[left] at (0,-3) {$4$};
\node[left] at (0,-4) {$5$};
\node[left] at (0,-5) {$6$};
\def\a{0.15};
\def\r{0.3};
\draw(0,0)--(1-\a,0)--(1+\a,-1)--(2-\r,-1);
\draw (2,-1-\r) arc (0:90:\r);
\draw[->](2,-1-\r)--(2,-5.5);
\draw(0,-1)--(1-\a,-1)--(1+\a,0)--(4-\a,0)--(4+\a,-1)--(6-\a,-1)--(6+\a,-2)--(7,-2);
\draw(0,-2)--(2-\a,-2)--(2+\a,-1)--(3-\a,-1)--(3+\a,-2)--(5-\r,-2);
\draw (5,-2-\r) arc (0:90:\r);
\draw[->](5,-2-\r)--(5,-5.5);
\draw(0,-3)--(2-\a,-3)--(2+\a,-2)--(3-\a,-2)--(3+\a,-1)--(4-\a,-1)--(4+\a,0)--(7,0);
\draw(0,-4)--(2-\a,-4)--(2+\a,-3)--(5-\a,-3)--(5+\a,-2)--(6-\a,-2)--(6+\a,-1)--(7,-1);
\draw(0,-5)--(2-\a,-5)--(2+\a,-4)--(5-\a,-4)--(5+\a,-3)--(7,-3);
\end{tikzpicture}
\caption{A reduced generalized wiring diagram}
\label{fig:generalized-wiring-01}
\end{figure}

Let $\mathcal{W}(k,n)$ be the set of reduced generalized wiring diagrams in which only wires labeled $1$ to $n$ ever occupy the first $k$ levels.
\begin{lem}
The maximum number of intersections on level $k$ among diagrams in $\mathcal{W}(k,n)$ equals $\Ans(k,n)$.
\end{lem}
\begin{proof}
Let $m(k,n)$ be the maximum number of intersections on level $k$ among $\mathcal{W}(k,n)$. Then we have $m(k,n)\geq\Ans(k,n)$ since a reduced word in $S_n$, which can be viewed as a wiring diagram, is an instance of $\mathcal{W}(k,n)$. 

Now for $W\in\mathcal{W}(k,n)$, at each timestamp $t$, let $A^{(t)}\subset{[n]\choose k}$ be the set of wires that occupy the first $k$ levels. Let $(A_0,A_1,\ldots,A_N)$ be the sequence of $A^{(t)}$'s where $A_i=A^{(t_i)}$ if step $t_i$ creates a crossing at level $k$. With the same reasoning as in Proposition~\ref{prop:equivalence-between-words-monotone-path}, we show that $(A_0,A_1,\ldots,A_N)$ is a monotone weakly separated path. Since $W$ is reduced, to go from $A_i$ to $A_{i+1}$, we take away some wire $a$ and add in some wire $b$ with $a<b$ and this shows monotonicity. Then for $i<j$, if there exists $b\in A_i\setminus A_j$, $a\in A_j\setminus A_i$ with $a<b$, then wire $a$ must cross wire $b$ (strictly) prior to step $t_i$ and also (strictly) between $t_i$ and $t_j$, contradicting $W$ being reduced. This shows that $(A_0,A_1,\ldots,A_N)$ is weakly separated. Here, $N$ is the number of crossings at level $k$ in $W$. By Proposition~\ref{prop:equivalence-between-words-monotone-path}, there exists a reduced word in $S_n$ with $N$ copies of $s_k$. Thus $m(k,n)\leq\Ans(k,n)$ as desired.
\end{proof}

We now consider a particularly nice set of reduced generalized wiring diagrams.
\begin{defn}\label{def:simple}
A reduced generalized wiring diagram $W$ is \textit{simple}, if there is no
\begin{itemize}
\item[(S1)] Pair of wires $a<b$ which intersects on a level other than $k$, and at the moment of their intersection wires $\{a+1,\ldots,b-1\}$ have already fallen; or
\item[(S2)] Wire $a$, which intersects with $k$ wires with larger labels (not counting when wire $a$ is going to infinity).
\end{itemize}
\end{defn}
We remark that (S1) implies that in a simple diagram, wire $a$ and wire $a+1$ can only intersect at level $k$. Also note that (S1) needs to be considered when the wire $a$ is falling to infinity. In particular, this means that in a simple diagram, the first step can only be an intersection between wires $k$ and $k+1$, or wire $k$ falling to infinity.

We write $\tilde{\mathcal{W}}(k,n)$ for the set of simple reduced generalized wiring diagrams from $\mathcal{W}(k,n)$. And write $\tilde{\mathcal{W}}(k)=\bigcup_{n>k} \tilde{\mathcal{W}}(k,n)$ for the set of simple reduced generalized wiring diagrams for fixed $k$.
\begin{proposition}\label{prop:simplify-diagram}
A diagram $W\in \mathcal{W}(k,n)$ can be transformed into a simple reduced diagram $W'\in\tilde{\mathcal{W}}(k,n)$ without changing the number of intersections on level $k$.
\end{proposition}
\begin{proof}
Let $W\in\mathcal{W}(k,n)$ and let $t$ be the first timestamp where $W$ violates some condition in Definition~\ref{def:simple}. If condition (S1) is violated by wires $a<b$ crossing normally (not during while $a$ is going to infinity), we simply remove this intersection to obtain $\tilde W$. The new diagram $\tilde{W}$ is still reduced, because in order for some wire $c$ to intersect the new wire $a$ (or $b$) twice without intersecting the previous $a$ and $b$ twice, $c$ must be between $a$ and $b$. However, (S1) says that $c$ has already fallen, so there are no such possibilities.

If condition (S1) is violated by wires $a<b$ intersecting as $a$ goes to infinity, we make sure that wire $b$ is the ``first" violation of condition (S1), i.e. the highest (with the smallest level). We expand this step of $a$ going to infinity by letting wire $a$ intersect those wires below $a$ and above $b$ at this timestamp first, and then going to infinity. Next, as above, we uncross the intersection between wires $a$ and $b$ by letting wire $b$ go to infinity instead while wire $a$ in $\tilde{W}$ takes on the role of wire $b$ in $W$ after this timestamp. The same argument in the last paragraph shows that $\tilde{W}$ is reduced. At the same time, the number of intersections at level $k$ stays unchanged.

If condition (S2) is violated with wire $a$, then we assume that at time $t$, wire $a$ crosses with $b_k$ at level $h$ where $a<b_k$; moreover, wire $a$ has already intersected with $b_1,\ldots,b_{k-1}$ that are larger than $a$. Let $\tilde{W}$ be obtained from $W$ by replacing the intersection at time $t$ with wire $a$ going to infinity. By reducedness, at time $t$, wires $b_1,\ldots,b_{k-1}$ must be at a higher (smaller) level than wire $a$, so $h\geq k$. This says that the number of intersections at level $k$ is the same in $\tilde W$. If $\tilde{W}$ becomes not reduced, then there must be some wire $c<a$ which is at level $\geq h+2$ at time $t$ so that it intersects $a$ the second time in $\tilde{W}$ at time $t$. However, since $c<b_1,\ldots,b_k$ and $c$ is at a lower (greater in value) level at time $t$, by condition (S2) and the minimality of $t$, wire $c$ must have fallen already. As a result, $\tilde W$ stays reduced.

We can continue the above process so that the end result $\tilde{W}$ is simple.
\end{proof}

For a diagram $W\in \mathcal{W}(k,n)$ and a timestamp $t$, we can associate a permutation $\pi^{(t)}_W=\pi^{(t)}\in S_{\infty}$ to it that records the positions of the non-fallen wires. To be precise, if $a_1<a_2<\cdots$ are the labels of the non-fallen wires at time $t$, then $a_{\pi^{(t)}(h)}$ is at level $h$ for $h=1,2,\ldots$. Here, the infinite symmetric group $S_{\infty}$ is the set of bijections on $\mathbb{Z}_{>0}$ with all but finitely many fixed points. We also let $f^{(t)}_W=f^{(t)}$ be the number of fallen wires of $W$ at timestamp $t$ and let $\kappa^{(t)}_W=\kappa^{(t)}$ be the number of intersections at level $k$ that have happened. In particular, we always start with $f^{(0)}=0$, $\kappa^{(0)}=0$ and $\pi^{(0)}=\id$.
\begin{lem}\label{lem:encoding-of-wiring-diagram}
A simple reduced generalized wiring diagram $W$ can be uniquely encoded by the sequence $\{(f_W^{(t)},\kappa_W^{(t)},\pi_W^{(t)})\}_{t}$ defined above. In other words, given a sequence $\{(f^{(t)},\kappa^{(t)},\pi^{(t)})\}_{t}$, there is at most one $W\in\tilde{\mathcal{W}}(k)$ such that $\{(f_W^{(t)},\kappa_W^{(t)},\pi_W^{(t)})\}_{t}=\{(f^{(t)},\kappa^{(t)},\pi^{(t)})\}_{t}$.
\end{lem}
\begin{proof}
Fix $\{(f^{(t)},\kappa^{(t)},\pi^{(t)})\}_{t}$ and we will recover $W\in\tilde{\mathcal{W}}(k)$ step by step. Note that there is a lot of redundancy in this encoding, as the information from $f$ and $\pi$ are almost sufficient. 

At step $t>0$, if $f^{(t)}=f^{(t-1)}$ meaning no wires fall, we simply apply a crossing at level $h$ if $\pi^{(t)}=\pi^{(t-1)}s_h$. The critical case is that $f^{(t)}=f^{(t-1)}+1$ meaning that a wire falls at this step. Note that from a permutation $\pi^{(t-1)}$, it is possible that deleting an entry (and flattening the permutation) will result in the same permutation as deleting another entry. For example, if $\pi^{(t-1)}=\id$, deleting any entry and flattening the values to $1,2,\ldots$ will result in $\pi^{(t)}=\id$. In such cases, to uniquely reconstruct a simple diagram $W$, the conditions in Definition~\ref{def:simple} become important.

Suppose that at time $t$, letting the wire at level $a$ go to infinity will result in the permutation $\pi^{(t)}$, i.e. deleting the entry at index $a$ of $\pi^{(t-1)}$ and flattening the values to $1,2,\ldots$ give us $\pi^{(t)}$, and letting the wire at level $b>a$ go to infinity will result in the same permutation $\pi^{(t)}$. Choose such minimal $a$ and maximal $b$. We analyze the permutation $u=\pi^{(t-1)}\in S_{\infty}$. 

First, for every positive integer $i$ such that $i<a$ or $i>b$, $u(i)$ must not lie in between $u(a)$ and $u(b)$. Secondly, if $u(a)>u(a+1)$, then letting this wire at level $a$ go to infinity results in a double crossing with the wire at level $a+1$ at this timestamp. Thus, $u(a)<u(a+1)$. By comparing the two permutations obtained from $u$ by deleting index $a$ and $b$ respectively, we must have $u(a+1)<u(a+2)$. This further implies $u(a+2)<u(a+3)$ and so on. Thus, $u(a)<u(a+1)<\cdots<u(b)$. It is now clear that deleting any index between $a$ and $b$ from $u$ results in the same permutation $\pi^{(t)}$, and finally we claim that at most one choice is possible. For $a\leq c\leq b-1$, if wire $c$ falls at time $t$, an intersection at level $c$ between this wire and the wire at level $c+1$ is created. By the arguments above, all the wires with labels between these two must have fallen (since they cannot exist before level $a$ or after level $b$), and by condition (S1) for simple diagrams, $c=k$. Thus, if $b\leq k$, only wire $b$ is allowed to fall; if $a\leq k<b$, only wire $b$ can fall if $\kappa^{(t)}=\kappa^{(t-1)}$ and only wire $k$ can fall if $\kappa^{(t)}=\kappa^{(t-1)}+1$; and the case $a>k$ cannot result in any valid diagrams.
\end{proof}

\begin{example}
We consider one optimal repeatable pattern of $k=2$, \[ 12 - 13 - 23 - 34 - 35 - 45 - 56 - 57 - 67 - 78 - 79 - 89 - \cdots, \]
discussed in Theorem~\ref{thm:M2n} and shown in Figure~\ref{fig:k=2-wiring}, and use a (simple) reduced generalized wiring diagram to describe it, shown in Figure~\ref{fig:k=2-generalized-wiring}.
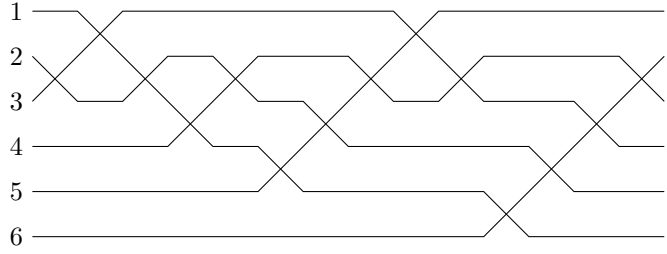
\begin{figure}[h!]
\centering
\begin{tikzpicture}[scale=0.6]
\node[left] at (0,0) {$1$};
\node[left] at (0,-1) {$2$};
\node[left] at (0,-2) {$3$};
\node[left] at (0,-3) {$4$};
\node[left] at (0,-4) {$5$};
\node[left] at (0,-5) {$6$};
\def\a{0.15};
\draw(0,0)--(1,0)--(2-\a,0)--(2+\a,-1)--(3-\a,-1)--(3+\a,-2)--(4-\a,-2)--(4+\a,-3)--(6-\a,-3)--(6+\a,-4)--(11-\a,-4)--(11+\a,-5)--(15,-5);
\draw(0,-1)--(1-\a,-1)--(1+\a,-2)--(3-\a,-2)--(3+\a,-1)--(5-\a,-1)--(5+\a,-2)--(7-\a,-2)--(7+\a,-3)--(12-\a,-3)--(12+\a,-4)--(15,-4);
\draw(0,-2)--(1-\a,-2)--(1+\a,-1)--(2-\a,-1)--(2+\a,0)--(9-\a,0)--(9+\a,-1)--(10-\a,-1)--(10+\a,-2)--(13-\a,-2)--(13+\a,-3)--(15,-3);
\draw(0,-3)--(4-\a,-3)--(4+\a,-2)--(5-\a,-2)--(5+\a,-1)--(8-\a,-1)--(8+\a,-2)--(10-\a,-2)--(10+\a,-1)--(14-\a,-1)--(14+\a,-2)--(15,-2);
\draw(0,-4)--(6-\a,-4)--(6+\a,-3)--(7-\a,-3)--(7+\a,-2)--(8-\a,-2)--(8+\a,-1)--(9-\a,-1)--(9+\a,0)--(15,0);
\draw(0,-5)--(11-\a,-5)--(11+\a,-4)--(12-\a,-4)--(12+\a,-3)--(13-\a,-3)--(13+\a,-2)--(14-\a,-2)--(14+\a,-1)--(15,-1);
\end{tikzpicture}
\caption{An optimal repeatable pattern for $k=2$}
\label{fig:k=2-wiring}
\end{figure}

In particular, the simple reduced generalized wiring diagram $\tilde{W}$ in Figure~\ref{fig:k=2-generalized-wiring} can be obtained from the wiring diagram $W$ in Figure~\ref{fig:k=2-wiring} via the simplification procedure in Proposition~\ref{prop:simplify-diagram}. Observe that the permutations $\pi^{(t)}_{\tilde{W}}$ are $\id,132,312,21,\id,132,312,21,\ldots$, which are periodic with period $4$. 

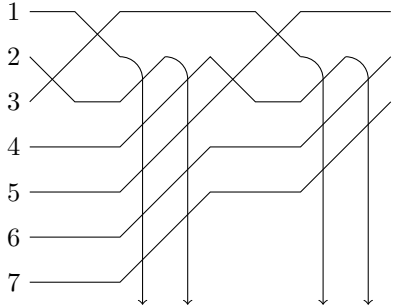
\begin{figure}[h!]
\centering
\begin{tikzpicture}[scale=0.6]
\def\a{0.15};
\def\r{0.3};
\node[left] at (0,0) {$1$};
\node[left] at (0,-1) {$2$};
\node[left] at (0,-2) {$3$};
\node[left] at (0,-3) {$4$};
\node[left] at (0,-4) {$5$};
\node[left] at (0,-5) {$6$};
\node[left] at (0,-6) {$7$};
\draw(0,0)--(2-\a,0)--(2+\a,-1)--(3-\r,-1);
\draw (3,-1-\r) arc (0:90:\r);
\draw[->](3,-1-\r)--(3,-6.5);
\draw(0,-1)--(1-\a,-1)--(1+\a,-2)--(3-\a,-2)--(3+\a,-1)--(4-\r,-1);
\draw (4,-1-\r) arc (0:90:\r);
\draw[->](4,-1-\r)--(4,-6.5);
\draw(0,-2)--(1-\a,-2)--(1+\a,-1)--(2-\a,-1)--(2+\a,0)--(6-\a,0)--(6+\a,-1)--(7-\r,-1);
\draw (7,-1-\r) arc (0:90:\r);
\draw[->](7,-1-\r)--(7,-6.5);
\draw(0,-3)--(3-\a,-3)--(3+\a,-2)--(4-\a,-2)--(4+\a,-1)--(5-\a,-1)--(5+\a,-2)--(7-\a,-2)--(7+\a,-1)--(8-\r,-1);
\draw (8,-1-\r) arc (0:90:\r);
\draw[->](8,-1-\r)--(8,-6.5);
\draw(0,-4)--(3-\a,-4)--(3+\a,-3)--(4-\a,-3)--(4+\a,-2)--(5-\a,-2)--(5+\a,-1)--(6-\a,-1)--(6+\a,0)--(9,0);
\draw(0,-5)--(3-\a,-5)--(3+\a,-4)--(4-\a,-4)--(4+\a,-3)--(7-\a,-3)--(7+\a,-2)--(8-\a,-2)--(8+\a,-1)--(9,-1);
\draw(0,-6)--(3-\a,-6)--(3+\a,-5)--(4-\a,-5)--(4+\a,-4)--(7-\a,-4)--(7+\a,-3)--(8-\a,-3)--(8+\a,-2)--(9,-2);
\end{tikzpicture}
\caption{An optimal repeatable pattern for $k=2$ via simple reduced generalized wiring diagrams corresponding to Figure~\ref{fig:k=2-wiring}}
\label{fig:k=2-generalized-wiring}
\end{figure}
\end{example}

\subsection{Finiteness of configurations and proof of the main theorem}
Let $\mathcal{T}_k$ be the set of all possible permutations $\pi^{(t)}_W$ at all timestamps $t$ across all $W \in \tilde{\mathcal{W}}(k)$. In this section, we will show that $\mathcal{T}_k$ is finite and resolve the rationality of $c_k$.

For a simple reduced generalized wiring diagram $W\in\tilde{\mathcal{W}}(k)$, let $N_t(W)$ be the set of non-fallen wires at timestamp $t$ and let $E_t(W)\subset N_t(w)$ be the set of non-fallen wires at timestamp $t$ that were on the first $k$ levels at or before timestamp $t$.

\begin{lem}\label{lem:no-larger-intersection}
For $W\in\tilde{\mathcal{W}}(k)$ and $t$, any wire $a\in N_t(W)\setminus E_t(W)$ does not intersect wires with larger labels at time $t$ or earlier.
\end{lem}
\begin{proof}
Assume the opposite and let $a\in N_t(W)\setminus E_t(W)$ intersect some wire $b>a$ at time $t$ or earlier. Note that the condition $a\in N_t(W)\setminus E_t(W)$ means that wire $a$ has not yet fallen but it has never been to the first $k$ levels. Choose $t$ to be minimal and then choose $b$ to be minimal among all such wires. The minimality of $t$ means that wire $a$ and wire $b$ intersect at time $t$ and let $h>k$ be the level of this intersection.

Each wire $i$ strictly between $a$ and $b$ must intersect $a$ or $b$ at some time before $t$. If such a wire $i$ intersects $a$ at time $t'<t$, then the minimality of $b$ is violated. So such a wire $i$ must not intersect $a$ and must intersect $b$, at some timestamp $t'<t$. If it does not intersect $b$ by falling, then $i\in N_{t'}(W)\setminus E_{t'}(W)$ since it has never been to level $k$ or above, and the minimality of $t$ is violated. As a result, all such wires $i$ have fallen at time $t$, and since $h>k$, condition (S1) is violated.
\end{proof}

\begin{lem}\label{lem:NF-bound}
For $W\in\tilde{\mathcal{W}}(k)$ and $t$, $|E_t(W)|\leq k^2+2k$.
\end{lem}
\begin{proof}
Assume the opposite that $|E_t(W)|\geq k^2+2k+1$. Let $a_1<a_2<\cdots$ be the set of wires that have not fallen at timestamp $t$ and let \[z=\max\{i\:|\: a_i\in E_t(W)\}.\]
Consider the set \[A=\{i\:|\: a_i\in E_t(W),a_i\text{ has intersected }a_z\}.\]
Since the wire $a_z$ has been to the first $k$ levels, there are at most $k-1$ wires from $a_1,\ldots,a_{z_1}$ that do not intersect $a_z$. This means $|A|\geq (z-1)-(k-1)\geq k^2+k+1$.

Since $|A|\geq k(k+1)+1$, the Erd\H{o}s--Szekeres Theorem says that we have either $k$ wires from $A$ which intersect pairwise, or $k+1$ wires from $A$ where no wires intersect. In the first case where a set $B\subset A$ of wires intersect pairwise with $|B|\geq k$, the smallest wire from $B$ then needs to intersect $k$ wires with larger labels, including $k-1$ wires from $B$ and $a_z$, contradicting (S2). In the second case where a set $B\subset A$ of wires have no intersections with $|B|\geq k+1$, the largest wire $a_i$ from $B$ can never visit the first $k$ levels, contradicting $a_i\in E_t(W)$. 
\end{proof}

\begin{corollary}\label{cor:T-finite}
For any $k$, $|\mathcal{T}_k|\leq k^{k^2+2k}$. In particular, $|\mathcal{T}_k|$ is finite.
\end{corollary}
\begin{proof}
Let $W\in\tilde{\mathcal{W}}(k)$ be simple and consider $\pi=\pi^{(t)}_W\in \mathcal{T}_k$. As above, let $a_1,a_2,\ldots$ be the labels of wires that have not yet fallen. Let $z=|E_t(W)|\leq k^2+2k$ and consider the Lehmer code $\code(\pi)$ where $\code(\pi)_i$ equals the number of wires $a_j$ that have intersected $a_i$, where $j>i$. It is a classical fact that Lehmer codes uniquely characterize permutations in $S_{\infty}$.

If $a_i\notin E_t(W)$, meaning that $a_i$ has not been to the first $k$ levels, then by Lemma~\ref{lem:no-larger-intersection}, any $a_j$ with $j>i$ does not intersect $a_i$ so has not been to the first $k$ levels either. This means that $E_t(W)=\{a_1,a_2,\ldots,a_z\}$ and that $\code(\pi)_i=0$ for $i>z$. At the same time, by (S2), each $a_i$ with $i\leq z$ can only intersect at most $k-1$ wires with larger labels. So $\code(\pi)_i\in\{0,1,\ldots,k-1\}$ for $i\in\{1,2,\ldots,z\}$. As a result, the total number of possible permutations is bounded by $k^z\leq k^{k^2+2k}$. 
\end{proof}

A \emph{piece} $P$ is a segment of a generalized wiring diagram $W$ containing the single move (adjacent crossing or falling wire) occurring in $W$ at some time $t$ together with the information of the permutations $\pi^-_P \coloneqq \pi_W^{(t-1)}$ and $\pi^+_P \coloneqq \pi_W^{(t)}$. The piece $P$ is \emph{simple} (with respect to $k$) if it can be obtained from a simple diagram $W \in \tilde{\mathcal{W}}(k)$. A series $P_1,\ldots,P_r$ of pieces such that $\pi^+_{P_i}=\pi^-_{P_{i+1}}$ for all $i$ may be concatenated into a \emph{pattern} $Q=P_1+\cdots+P_r$, which is the segment of a generalized wiring diagram obtained by drawing $P_1,\ldots,P_r$ next to each other, together with the information of $\pi_Q^-\coloneqq\pi_{P_1}^-$ and $\pi_Q^+ \coloneqq\pi_{P_r}^+$. 

The following proposition shows that being simple is a local property of a generalized wiring diagram: a diagram is simple if and only if all of its constituent pieces are simple.

\begin{prop}
\label{prop:simple-pattern-is-simple}
Let $Q=P_1 + \cdots + P_r$, where each piece $P_i$ is simple with respect to $k$ and where $\pi_{P_1}^-=\id$, then $Q \in \tilde{\mathcal{W}}(k)$. 
\end{prop}
\begin{proof}
First note that $Q$, viewed as a generalized wiring diagram, is reduced, since each piece, by virtue of coming from a reduced diagram and carrying with it the permutations $\pi^-$ and $\pi^+$ clearly preserves reducedness when concatenated. The condition (S1) from Definition~\ref{def:simple} is also clearly preserved when we apply each piece, since it is equivalent to the condition that no simple piece swaps two wires with adjacent labels, or has a wire fall from a level other than $k$ when the wire below it has label one larger, where we read labels from $\pi^-$. Finally, to check condition (S2), we can just check that for each index $i$ and time $t$, 
\[
|\{j<i \mid \pi_Q^{(t)}(j)>\pi_Q^{(t)}(i) \}| < k.
\]
This is because, since $Q$ is reduced, none of the wires with label higher than $a_{\pi_Q^{(t)}(i)}$ which have crossed this wire can have fallen at or before time $t$, so $|\{j<i \mid \pi_Q^{(t)}(j)>\pi_Q^{(t)}(i) \}|$ is this number of wires. Since this is a condition satisfied by all $\pi^{(t)}_W$ for $W$ simple, it is satisfied by $\pi_Q^{(t)}=\pi_{P_{t+1}}^-$.
\end{proof}

We call a pattern $Q=P_1 + \cdots +P_r$ \emph{simple} (with respect to $k$) if all of its constituent pieces $P_i$ are simple with respect to $k$. By Proposition~\ref{prop:simple-pattern-is-simple}, this does not conflict with our earlier definition of simple diagrams.

\begin{defn}
\label{def:K}
For $f \in \mathbb{N}$ and $\pi \in \mathcal{T}_k$, let $K(f,\pi)$ be the maximum number of crossings at level $k$ among all simple patterns $Q=P_1+\cdots+P_r$ such that $\pi^-_{P_1},\ldots,\pi^-_{P_r}$ are distinct, $\pi^-_{P_1}=\pi^+_{P_r}=\pi$, and $Q$ has $f$ fallen wires. Since $\mathcal{T}_k$ is finite by Corollary~\ref{cor:T-finite}, there are finitely many simple patterns whose constituent pieces have distinct values of $\pi^-$.
\end{defn}

\begin{thm}
\label{thm:ck-is-rational}
For any $k \in \mathbb{Z}_{>0}$ we have:
\begin{equation}
\label{eq:ck-equals-max}
c_k = \max_{f \in \mathbb{N}, \pi \in \mathcal{T}_k} \frac{K(f,\pi)}{f}.
\end{equation}
In particular, $c_k$ is rational.
\end{thm}
\begin{proof}
We first show that $c_k \geq \max_{f \in \mathbb{N}, \pi \in \mathcal{T}_k} \frac{K(f,\pi)}{f}$. Let $(f_0,\pi_0)$ be such that $K(f_0,\pi_0)/f_0$ achieves the maximum, and let $Q$ be a simple pattern realizing $K(f_0,\pi_0)$ crossings at level $k$, with $\pi_Q^-=\pi_Q^+$, as in Definition~\ref{def:K}. Choose $R \in \tilde{\mathcal{W}}(k)$ with $\pi_R^-=\id, \pi_R^+=\pi_Q^-,$ and no fallen wires. For any $m \in \mathbb{N}$, consider the diagram
\[
D_m \coloneqq R + mQ = R + \underbrace{Q + \cdots + Q}_{\text{$m$ copies}}.
\]
By Proposition~\ref{prop:simple-pattern-is-simple}, we have $D_m \in \tilde{\mathcal{W}}(k)$ since $R$ and $Q$ are simple. By Lemma~\ref{lem:NF-bound}, and since $Q$ has $f_0$ fallen wires, we have $D_m \in \tilde{\mathcal{W}}(k,k^2+2k+mf_0)$. The number of intersections on level $k$ of $D_m$ is the number $r$ in $R$, plus $mK(f_0,\pi_0)$. We conclude
\[
c_k \geq \lim_{m \to \infty} \frac{r+mK(f_0,\pi_0)}{k^2+2k+mf_0} =  \frac{K(f_0,\pi_0)}{f_0}=\max_{f \in \mathbb{N}, \pi \in \mathcal{T}_k} \frac{K(f,\pi)}{f}.
\]

We now prove the upper bound on $c_k$. Write $M$ for the maximum on the right-hand side of (\ref{eq:ck-equals-max}), and let $W$ be any diagram from $\tilde{\mathcal{W}}(k)$. Express $W$ uniquely as a sum of its constituent simple pieces:
\[
W=P_1+\cdots+P_r.
\]
Note that $W$ has $f^{(r)}=f^{(r)}_W$ fallen wires and $\kappa^{(r)}=\kappa^{(r)}_W$ crossings at level $k$.
Let $i \leq j$ be a closest pair of indices such that $\pi^-_{P_i}=\pi^+_{P_j}$, and write $\pi$ for this common permutation. Consider the encodings $(f^{(i-1)},\kappa^{(i-1)},\pi)$ and $(f^{(j)},\kappa^{(j)},\pi)$ immediately before and after the pieces $P_i,P_j$ respectively. Since $\pi^+_{P_{i-1}}=\pi^-_{P_{j+1}}$, we may form a new diagram $W'=P_1+\cdots+P_{i-1}+P_{j+1}+\cdots+P_r$ by removing the piece $P_i+\cdots+P_j$ and $W' \in \tilde{\mathcal{W}}(k)$ by Proposition~\ref{prop:simple-pattern-is-simple}. This new diagram has $f^{(r)}-(f^{(j)}-f^{(i-1)})$ fallen wires and $\kappa^{(r)}-(\kappa^{(j)}-\kappa^{(i-1)})$ crossings at level $k$. By construction, we have
\[
\kappa^{(j)}-\kappa^{(i-1)} \leq (f^{(j)}-f^{(i-1)})M.
\]
Since $W'$ remains simple, we can iteratively apply this procedure until we reach a diagram $W_0$ with all permutations distinct. By Corollary~\ref{cor:T-finite}, there are finitely many such diagrams, so there is some maximum possible number $K_0$ of level-$k$ crossings in $W_0$. Since $W \in \mathcal{W}(k,k^2+2k+f^{(r)})$, we see that
\[
\mathcal{M}(k,k^2+2k+f^{(r)}) \leq K_0 + f^{(r)}M.
\]
Dividing by $k^2+2k+f^{(r)}$ and letting $f^{(r)} \to \infty$, we obtain
\[
c_k \leq M,
\]
as desired.
\end{proof}

\begin{corollary}
\label{cor:ck-achieved-by-repeatable}
The optimal density $c_k$ can be achieved by a repeatable pattern, in the sense of Section~\ref{sec:bounds}.
\end{corollary}
\begin{proof}
Let $(f_0,\pi_0)$ achieve the maximum in the right-hand side of (\ref{eq:ck-equals-max}), and let $Q$ be the diagram realizing $f_0,\pi_0$ as in the proof of Theorem~\ref{thm:ck-is-rational}, which, by the theorem, achieves the density $c_k$. Taking the labels of the wires occupying the top $k$ levels in a reduced generalized wiring diagram gives a monotone weakly separated path by the same reasoning as in Proposition~\ref{prop:equivalence-between-words-monotone-path}, and applying this operation to $Q$ gives the desired repeatable pattern, by construction.
\end{proof}

\begin{prop}
    Let $p_k(n) \coloneqq \mathcal{M}(k,n)-c_kn$. Then, for sufficiently large $n$, $p_k(n)$ is periodic in $n$.
\end{prop}
\begin{proof}
    Let $\mathcal{M}(k, n, \pi)$ be the maximal number of intersections on level $k$ among diagrams from $\tilde{W}(k,n)$ having final timestamp $\pi$. Clearly $\mathcal{M}(k, n) = \max_{\pi \in \mathcal{T}_k} \mathcal{M}(k, n, \pi)$. Define $p_k(n, \pi) \coloneqq M(k, n, \pi) - c_kn$.

    It follows from the proof of Theorem~\ref{thm:ck-is-rational} that $p_k(n,\pi)$ is a bounded function of $n$ and $\pi$. Furthermore, since $c_k$ is rational, this function takes rational values of bounded denominator, and thus $p_k(n,\pi)$ attains only finitely many different values. Thus we may find $a<b \in \mathbb{N}$ so that $p_k(a,\pi)=p_k(b,\pi)$ for all $\pi \in \mathcal{T}_k$, since $\mathcal{T}_k$ is finite by Corollary~\ref{cor:T-finite}. Now, it is clear that $p_k(n,\pi)$ depends only on $\{p_k(n-1,\sigma)\}_{\sigma \in \mathcal{T}_k}$. Thus we have that $p_k(n+b-a,\pi)=p_k(n,\pi)$ for all $n \geq a$ and all $\pi \in \mathcal{T}_k$. In particular, for $n$ sufficiently large, $p_k(n)=\max_{\pi \in \mathcal{T}_k} p_k(n,\pi)$ is periodic in $n$.
\end{proof}

\section{The minimization problem for finite Coxeter groups}
\label{sec:other-types}
In this section, we investigate a related question: for the longest element $w_0$ of a finite Coxeter group $W$, what is the minimum number of appearances of a generator $s_i$ in $\Red(w_0)$, the set of reduced words for $w_0$. This question is very easy in type $A_{n-1}$ where $W\simeq\mathfrak{S}_n$. Namely, the minimum number of occurrences of the simple transposition $(i\ i+1)$ in $\Red(w_0)$ is $\min\{i, n-i\}$. We observe a surprising phenomenon with respect to these numbers and the Cartan matrix of $W$ (Theorem~\ref{thm:Vinberg-main}).

Throughout this section, let $$W=\langle s_1,\ldots,s_n\:|\: (s_is_j)^{m_{ij}}=\mathrm{id}\text{ for all }i,j\rangle$$ be a finite Coxeter group generated by a set of simple reflections $S=\{s_1,\ldots,s_n\}$. For $w\in W$, let $\ell(W)$ denote the Coxeter length of $w$. For $J \subseteq S$, the \emph{parabolic subgroup} $W_J$ is the subgroup of $W$ generated by $J$, viewed as a Coxeter group with simple reflections $J$.  Each left coset $wW_J$ of $W_J$ in $W$ contains a unique element $w^J$ of minimal length, and the set $\{w^J \: | \: w \in W\}$ of these minimal coset representatives is called the \emph{parabolic quotient} $W^J$.  Letting $w_J \in W_J$ be the unique element such that $w^Jw_J=w$, we have $\ell(w^J)+\ell(w_J)=\ell(w)$ and this is called the \emph{parabolic decomposition} of $w$. As $W$ is finite, $W^J$ is finite and it contains a unique element $w_0^J$ of maximum length. We utilize the Bruhat order on $W$ and $W^J$, where $u\leq w$ if $u$ equals a subword of a (or equivalently, any) reduced word of $w$. For convenience, we adopt the notation that $J_i:=S-\{s_i\}$ for each $s_i\in S$. We refer readers to \cite{bjorner-brenti} for a detailed exposition on Coxeter groups.

We start with an algorithm to compute the minimum number of $s_i$ that appears in $\Red(w)$ for all $w$. 
\begin{prop}\label{prop:min-algorithm}
Fix $w\in W$ and $s_i\in S$. Define a sequence of Coxeter group elements $w^{(0)},w^{(1)},\ldots$ as follows: $w^{(0)}=w^{J_i}$ and $w^{(k+1)}=(w^{(k)}s_i)^{J_i}$ if $w^{(k)}\neq\mathrm{id}$, for $k\geq0$. This algorithm will eventually stop (at some $w^{(N)}=\mathrm{id}$). Then the minimum number of $s_i$ that appears in $\Red(w)$ is the $k$ for which $w^{(k)}=\mathrm{id}$.
\end{prop}
\begin{proof}
First notice that in this procedure, if $w^{(j)}\neq\mathrm{id}$, then as $w^{(j)}\in W^{J_i}$, it must have a single descent at $s_i$. As a result, $\ell(w^{(j+1)})\leq\ell(w^{(j)}s_i)<\ell(w^{(j)})$ so we will eventually end up at the identity. This procedure also produces a (class of) reduced word of $w$ with $k$ $s_i$'s where $w^{(k)}=\mathrm{id}$.

Let $k$ be such that $w^{(k)}=\mathrm{id}$ and take an arbitrary reduced word $s_{i_1}s_{i_2}\cdots s_{i_{\ell}}$ of $w$. Pick out the $s_i$'s in this reduced word as $i_{a_K}=i_{a_{K-1}}=\cdots=i_{a_1}=i$ where $a_K<a_{K-1}<\cdots<a_1$. For $j=0,1,\ldots,K-1$, let $u^{(j)}=s_{i_1}s_{i_2}\cdots s_{i_{a_{j+1}}}$ which is the product from $s_{i_1}$ to the $(j+1)^{th}$ $s_i$ in this reduced word counted from the right. Also say $u^{(K)}=\mathrm{id}$. 

Recall the following standard fact of Coxeter groups: if $x\leq y$, then $x^J\leq y^J$ for any subset $J\subset S$. This can be proved via an application of the subword property of Bruhat orders. Also see \cite{bjorner-brenti}.

We now show that $u^{(j)}\geq w^{(j)}$ for $j=0,1,\ldots,k$ in the Bruhat order by induction. For the base case, notice that both $u^{(0)}$ and $w^{(0)}$ are in the left coset $wW_{J_i}$ and since $w^{(0)}$ is the minimal coset representative, we have $u^{(0)}\geq w^{(0)}$. Now assume $u^{(j)}\geq w^{(j)}\neq\mathrm{id}$ for some $j\geq0$. By definition, both of them have a right descent at $s_i$ so we have $u^{(j)}s_i\geq w^{(j)}s_i$ by the fact in the last paragraph with $J=\{s_i\}$. With another application of this fact with $J=J_i$, we have $(u^{(j)}s_i)^{J_i}\geq (w^{(j)}s_i)^{J_i}=w^{(j+1)}$. At the same time, $u^{(j+1)}$ and $u^{(j)}s_i$ are in the same coset of $W_{J_i}$ by definition, so $u^{(j+1)}\geq (u^{(j)}s_i)^{J_i}\geq w^{(j+1)}$. The induction step goes through.

Finally, $u^{(k-1)}\geq w^{(k-1)}\neq\mathrm{id}$. This means $u^{(k-1)}\neq\mathrm{id}$ so $K>k-1$, $K\geq k$ as desired.
\end{proof}

Recall that a \emph{generalized Cartan matrix} $A$ of a Coxeter system $(W,S)$ is a real $n\times n$ matrix such that
\begin{itemize}
\item $A_{ii}=2$ for $i=1,\ldots,n$ and $A_{ij}\leq0$ for $i\neq j$,
\item $A_{ij}<0$ if and only if $A_{ji}<0$ and $A_{ij}A_{ji}=4\cos^2(\pi/m_{ij})$ for $i\neq j$.
\end{itemize}
We say that a generalized Cartan matrix $A$ is \emph{restricted} if $m_{ij}=3$, or equivalently, there is a single edge between $s_i$ and $s_j$ in the Dynkin diagram, implies that $A_{ij}=A_{ji}=-1$. Note that if $(W,S)$ is simply-laced, then any restricted generalized Cartan matrix is the Cartan matrix. We now state our main result of the section.

\begin{thm}\label{thm:Vinberg-main}
Let $W$ be a finite Coxeter group generated by $S=\{s_1,\ldots,s_n\}$ and let $v\in\mathbb{R}_{>0}^n$ be such that $v_i$ is the minimum number of appearances of $s_i$ in a reduced word of $w_0$. Then there exists a restricted generalized Cartan matrix $A\in\mathbb{R}^{n\times n}$ of $(W,S)$ such that $Av\geq\mathbf{0}$, where the comparison is made entry-wise.
\end{thm}
\begin{proof}
We make use of Proposition~\ref{prop:min-algorithm} for each type separately and provide the corresponding restricted generalized Cartan matrix. Note that $(Av)_i=2v_i+\sum_{j\sim i}A_{ij}v_j$, where $j\sim i$ means that the nodes $i$ and $j$ are adjacent in the Dynkin diagram. So $Av\geq\mathbf{0}$ is intuitively saying that the value $v$ at each node $i$ is at least half of the weighted sum of its neighbors.

For the classical types, we mainly argue about type $B_n$, whose Coxeter group $W(B_n)$ is isomorphic to the group of signed permutations. The argument for type $D_n$, whose Coxeter group $W(D_n)$ is an index-2 subgroup of $W(B_n)$, is similar and we will omit unnecessary details. The argument for type $A_n$ is simpler. And for the exceptional types, we use Proposition~\ref{prop:min-algorithm} and a computer to generate each $v_i$, and then provide the matrix $A$ directly.
\vskip0.5em
\noindent\textbf{Type $B_n$}. Writing $\bar i:=-i$, we adopt the convention that 
\[W(B_n)=\{w\text{ is a permutation on }\bar n,\ldots,\bar1,1,2,\ldots,n\:|\: w(i)=-w(\bar i)\ \forall i\}\]
which is generated by $S=\{s_1,\ldots,s_n\}$ where $s_1=(1\ \bar1)$, $s_i=(i{-}1\ i)(\overline{i{-}1}\ \bar i)$ in cycle notation for $i=2,\ldots,n$. The Dynkin diagram and the $v_i$'s that we are about to compute can be seen in Figure~\ref{fig:minimal-Bn}.
\begin{figure}[h!]
\centering
\begin{tikzpicture}[scale=1.0]
\node at (0,0) {$\bullet$};
\node at (1,0) {$\bullet$};
\node at (2,0) {$\bullet$};
\node at (3,0) {$\bullet$};
\node at (4.5,0) {$\bullet$};
\node at (5.5,0) {$\bullet$};
\node at (6.5,0) {$\bullet$};
\node at (7.5,0) {$\bullet$};
\node at (3.75,0) {$\cdots$};
\draw(1,0)--(3,0);
\draw(4.5,0)--(7.5,0);
\draw(0,0.03)--(1,0.03);
\draw(0,-0.03)--(1,-0.03);
\node[above] at (0,0) {$s_1$};
\node[above] at (1,0) {$s_2$};
\node[above] at (2,0) {$s_3$};
\node[above] at (5.5,0) {$s_{n-2}$};
\node[above] at (6.5,0) {$s_{n-1}$};
\node[above] at (7.5,0) {$s_{n}$};
\node[above] at (-1,0) {labels};
\node[below] at (-1,0) {$v_i$'s};
\node[below] at (0,0) {$n$};
\node[below] at (1,0) {$n$};
\node[below] at (2,0) {$n{-}1$};
\node[below] at (3,0) {$n{-}2$};
\node[below] at (5.5,0) {$4$};
\node[below] at (6.5,0) {$3$};
\node[below] at (7.5,0) {$2$};
\end{tikzpicture}
\caption{minimal number of occurrences of each $s_i$ in reduced words of $w_0$ of type $B_n$, with labels on top and $v_i$'s on the bottom}
\label{fig:minimal-Bn}
\end{figure}
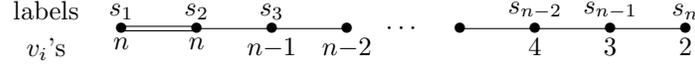

We write element $w\in W(B_n)$ in one-line notation given by $w(1) w(2)\cdots w(n)$. The longest element is $w_0=\bar1 \bar2 \cdots \bar n$ while the identity is $\mathrm{id}= 12\cdots n$. Fix some $s_i\in S$ with $i\geq2$ and we now run through the algorithm in Proposition~\ref{prop:min-algorithm}. Keep notations as in Proposition~\ref{prop:min-algorithm}, we use induction on $k$ to show that
\[w^{(k)}=1\ 2\ \cdots\  i{-}2\ n{+}1{-}k\  \overline{n{-}k}\ \overline{n{-}k{-}1}\ \cdots\ \bar i\ i{-}1\ n{-}k{+2}\ \cdots\ n\]
for $k\geq1$, where $\cdots$ indicates a sequence of consecutive increasing integers. We start with $w^{(0)}=12\cdots i{-}1\ \bar n\cdots \overline{i{+}1}\ \overline{i}$ and then \[w^{(0)}s_i=1\ 2\ \cdots\ i{-}2\ \bar n\ i{-}1\ \overline{n{-}1}\ \cdots\ \overline{i+1}\ \overline{i}.\]
Taking the parabolic quotient to obtain $w^{(1)}=(w^{(0)}s_i)^{J_i}$, where $J_i=S\setminus\{s_i\}$, we effectively get rid of the signs in coordinates $1,\ldots,i-1$ and sort these values, and also sort the values in coordinates $i,i+1,\ldots,n$ respectively. This gives
\[w^{(1)}=1\ 2\ \cdots\ i{-}2\ n\ \overline{n{-}1}\ \cdots\ \overline{i{+}1}\ \overline{i}\ i{-}1\]
as desired, establishing the base case. Checking the inductive steps is also done in the same way, by writing down
\[w^{(k)}s_i=1\ 2\ \cdots\  i{-}2\ \overline{n{-}k}\ n{+}1{-}k\ \overline{n{-}k{-}1}\ \cdots\ \bar i\ i{-}1\ n{-}k{+2}\ \cdots\ n,\]
sorting the values without the signs in coordinates $1,\ldots,i-1$ and sorting the values while keeping the signs in coordinates $i,\ldots,n$ to obtain $w^{(k+1)}=(w^{(k)}s_i)^{J_i}$. Finally, when $k=n+1-i$, we see that $w^{(k)}=s_i$ so we conclude that $v_i=n+2-i$. 

To see that $v_1=n$, we notice in fact that every time $s_1$ is applied to reduce the length of a signed permutation, there is one less negative values among $w(1),\ldots,w(n)$, and every time some other $s_i$ is applied, where $i\geq2$, the number of negative values among $w(1),\ldots,w(n)$ stays the same. This directly gives $v_1=n$. 

Finally, to specify a restricted generalized Cartan matrix $A\in\mathbb{R}^{n\times n}$, it suffices to specify $A_{12}=-2$ and $A_{21}=-1$. We check that $Av=(0,1,0,\ldots,0,0,1)^T\geq\mathbf{0}$. 

\vskip0.5em
\noindent\textbf{Type $D_n$}. The same argument as in type $B_n$ works in this case, by explicitly writing down the signed permutations $w^{(0)},w^{(1)},\ldots$ for each $s_i$. We omit the tedious details here and provide the answers for $v_i$'s in Figure~\ref{fig:minimal-Dn}.
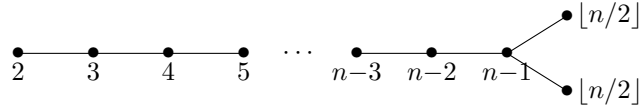
\begin{figure}[h!]
\centering
\begin{tikzpicture}
\node at (0,0) {$\bullet$};
\node at (1,0) {$\bullet$};
\node at (2,0) {$\bullet$};
\node at (3,0) {$\bullet$};
\draw(0,0)--(3,0);
\node at (3.75,0) {$\cdots$};
\node at (4.5,0) {$\bullet$};
\node at (5.5,0) {$\bullet$};
\node at (6.5,0) {$\bullet$};
\draw(4.5,0)--(6.5,0);
\node at (7.3,0.5) {$\bullet$};
\node at (7.3,-0.5) {$\bullet$};
\draw(7.3,0.5)--(6.5,0)--(7.3,-0.5);
\node[below] at (0,0) {$2$};
\node[below] at (1,0) {$3$};
\node[below] at (2,0) {$4$};
\node[below] at (3,0) {$5$};
\node[below] at (4.5,0) {$n{-}3$};
\node[below] at (5.5,0) {$n{-}2$};
\node[below] at (6.5,0) {$n{-}1$};
\node[right] at (7.3,0.5) {$\lfloor n/2\rfloor$};
\node[right] at (7.3,-0.5) {$\lfloor n/2\rfloor$};
\end{tikzpicture}
\caption{minimal number of occurrences of each $s_i$ in reduced words of $w_0$ of type $D_n$}
\label{fig:minimal-Dn}
\end{figure}
Since type $D_n$ is simply-laced, the restricted generalized Cartan matrix is fixed. We check that $Av\geq\mathbf{0}$, which in fact has value $0$ at most coordinates.

\vskip0.5em
\noindent\textbf{Type $A_{n-1}$}. We quickly go over the algorithm in Proposition~\ref{prop:min-algorithm}. Fix $i\leq \frac{n-1}{2}$. Let $Z_k=\{w^{(k)}(1),\ldots,w^{(k)}(i)\}$ so that $Z_0=\{n-i+1,\ldots,n\}$. To obtain $Z_{k+1}$ from $Z_k$, we exchange the largest entry of $Z_{k}$ with the smallest entry of $\{1,\ldots,n\}\setminus Z_{k}$. It is then immediate that $Z_i$ becomes $\{1,\ldots,i\}$ so that $v_i=i$. By symmetry of the Dynkin diagram, $v_i=n-i$ for $i\geq\frac{n-1}{2}$. We check again that most entries of $Av$ are zeroes, except one or two positive integers in the middle.

\vskip0.5em
\noindent\textbf{Type $E_6$, $E_7$, $E_8$.} The $v_i$'s are shown in Figure~\ref{fig:minimal-En}. We check that each $v_i$ is at least half of the sum of its neighbors.
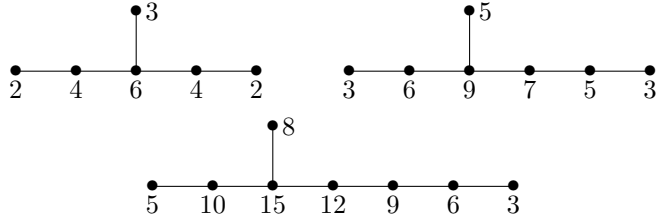
\begin{figure}[h!]
\centering
\begin{tikzpicture}[scale=0.8]
\node at (0,0) {$\bullet$};
\node at (1,0) {$\bullet$};
\node at (2,0) {$\bullet$};
\node at (3,0) {$\bullet$};
\node at (4,0) {$\bullet$};
\node at (2,1) {$\bullet$};
\draw(0,0)--(4,0);
\draw(2,0)--(2,1);
\node[below] at (0,0) {$2$};
\node[below] at (1,0) {$4$};
\node[below] at (2,0) {$6$};
\node[below] at (3,0) {$4$};
\node[below] at (4,0) {$2$};
\node[right] at (2,1) {$3$};
\end{tikzpicture}
\qquad
\begin{tikzpicture}[scale=0.8]
\node at (0,0) {$\bullet$};
\node at (1,0) {$\bullet$};
\node at (2,0) {$\bullet$};
\node at (3,0) {$\bullet$};
\node at (4,0) {$\bullet$};
\node at (5,0) {$\bullet$};
\node at (2,1) {$\bullet$};
\draw(0,0)--(5,0);
\draw(2,0)--(2,1);
\node[below] at (0,0) {$3$};
\node[below] at (1,0) {$6$};
\node[below] at (2,0) {$9$};
\node[below] at (3,0) {$7$};
\node[below] at (4,0) {$5$};
\node[below] at (5,0) {$3$};
\node[right] at (2,1) {$5$};
\end{tikzpicture}

\begin{tikzpicture}[scale=0.8]
\node at (0,0) {$\bullet$};
\node at (1,0) {$\bullet$};
\node at (2,0) {$\bullet$};
\node at (3,0) {$\bullet$};
\node at (4,0) {$\bullet$};
\node at (5,0) {$\bullet$};
\node at (6,0) {$\bullet$};
\node at (2,1) {$\bullet$};
\draw(0,0)--(6,0);
\draw(2,0)--(2,1);
\node[below] at (0,0) {$5$};
\node[below] at (1,0) {$10$};
\node[below] at (2,0) {$15$};
\node[below] at (3,0) {$12$};
\node[below] at (4,0) {$9$};
\node[below] at (5,0) {$6$};
\node[below] at (6,0) {$3$};
\node[right] at (2,1) {$8$};
\end{tikzpicture}
\caption{minimal number of occurrences of each $s_i$ in reduced words of $w_0$ of type $E_6$, $E_7$, $E_8$.}
\label{fig:minimal-En}
\end{figure}

\vskip0.5em
\noindent\textbf{Type $F_4$}. The $v_i$'s are shown in Figure~\ref{fig:minimal-F4} and we specify $A_{2,3}=A_{3,2}=-\sqrt2$.
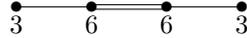
\begin{figure}[h!]
\centering
\begin{tikzpicture}
\node at (0,0) {$\bullet$};
\node at (1,0) {$\bullet$};
\node at (2,0) {$\bullet$};
\node at (3,0) {$\bullet$};
\draw(0,0)--(1,0);
\draw(2,0)--(3,0);
\draw(1,-0.03)--(2,-0.03);
\draw(1,0.03)--(2,0.03);
\node[below] at (0,0) {$3$};
\node[below] at (1,0) {$6$};
\node[below] at (2,0) {$6$};
\node[below] at (3,0) {$3$};
\end{tikzpicture}
\caption{minimal number of occurrences of each $s_i$ in reduced words of $w_0$ of type $F_4$}
\label{fig:minimal-F4}
\end{figure}

\noindent\textbf{Type $H_3$ and $H_4$}. The $v_i$'s are shown in Figure~\ref{fig:minimal-H} and we specify $A_{2,3}=-2\cos^2(\pi/5)$, $A_{3,2}=-2$ for type $H_3$, and $A_{3,4}=-2\cos^2(\pi/5)$, $A_{4,3}=-2$ for type $H_4$.
\begin{figure}[h!]
\centering
\begin{tikzpicture}
\node at (0,0) {$\bullet$};
\node at (1,0) {$\bullet$};
\node at (2,0) {$\bullet$};
\draw(0,0)--(2,0);
\node at (1.5,0) {$5$};
\node[below] at (0,0) {$3$};
\node[below] at (1,0) {$5$};
\node[below] at (2,0) {$5$};
\end{tikzpicture}
\qquad
\begin{tikzpicture}
\node at (0,0) {$\bullet$};
\node at (1,0) {$\bullet$};
\node at (2,0) {$\bullet$};
\node at (3,0) {$\bullet$};
\draw(0,0)--(3,0);
\node at (2.5,0) {$5$};
\node[below] at (0,0) {$5$};
\node[below] at (1,0) {$10$};
\node[below] at (2,0) {$15$};
\node[below] at (3,0) {$15$};
\end{tikzpicture}
\caption{minimal number of occurrences of each $s_i$ in reduced words of $w_0$ of type $H_3$ and $H_4$}
\label{fig:minimal-H}
\end{figure}

\noindent\textbf{Type $I_n$}. Here, $v_1=v_2=\lceil m_{12}/2\rceil$ so we let $A_{1,2}=-2\cos^2(\pi/m_{12})$ and $A_{2,1}=-2$.
\end{proof}

\section*{Acknowledgements}
We are grateful to the participants of the Harvard--MIT Combinatorics Preseminar, where some initial discussions of this problem took place, and in particular to Darij Grinberg for sharing this problem. We would also like to thank Sorawee Porncharoenwase for producing useful computer code for this project and the anonymous referees for their careful reading.

\bibliographystyle{plain}
\bibliography{arxiv-final.bib}
\end{document}